%% file: complexity.tex
\documentclass[11pt]{amsart}

\usepackage[dvips]{graphicx}
\usepackage[letterpaper,hmargin=3cm,vmargin=3cm,dvips]{geometry}
\usepackage{url}
\usepackage{mdwlist}
\usepackage{psfrag}
\usepackage{hyphenat}
\usepackage{booktabs}
\usepackage{mdwtab}

\input{macros}
\input{defs}
\begin{document}
\title{3-manifolds efficiently bound 4-manifolds}

\author[Costantino]{Francesco Costantino}
\thanks{FC was supported by a Marie Curie Fellowship issued by the
  European Community and hosted by Institut de Recherche
  Math\'emathique Avanc\'ee de Strasbourg.}
\address{Institute de Recherche Math\'ematique Avanc\'ee (IRMA),
  7, Rue Ren\'e Descartes\\
  67084 Strasbourg Cedex, France}
\email{costanti@math.u-strasbg.fr}

\author[Thurston]{Dylan Thurston}
\thanks{DPT~was supported by NSF Grant DMS-0071550, Harvard
University, the University of Pisa, and a Sloan Research Fellowship.}
\address {Department of Mathematics, Barnard College, Columbia University\\ New York, NY 10027}
\email {dthurston@barnard.edu}

\begin{abstract}
  It is known since 1954 that every 3-manifold bounds a 4-manifold.
  Thus, for instance, every 3-manifold has a surgery diagram.  There
  are several proofs of this fact, but
  there has been little attention to the complexity of the 4-manifold
  produced.  Given a 3-manifold $M^3$ of complexity $n$, we
  construct a 4-manifold bounded by $M$ of complexity $O(n^2)$, where
  the ``complexity'' of a piecewise-linear manifold is the
  minimum number of $n$-simplices in a triangulation.

  The proof goes through the notion of ``shadow complexity'' of a
  3-manifold~$M$.  A shadow of~$M$ is a well-behaved 2-dimensional
  spine of a 4-manifold bounded by~$M$.  We further prove that, for a
  manifold~$M$ satisfying the Geometrization Conjecture with Gromov
  norm~$G$ and shadow complexity~$S$,
  \[
  c_1G \leq S \leq c_2G^2
  \]
  for suitable constants $c_1$, $c_2$.  In particular, the manifolds
  with shadow complexity~0 are the graph manifolds.

  In addition, we give an $O(n^4)$ bound for the complexity of
  a spin 4-manifold bounding a given spin 3-manifold.  We also show
  that every stable map from a 3-manifold~$M$ with Gromov norm $G$ to
  $\RR^2$ has at least $G/10$ crossing singularities, and if $M$ is
  hyperbolic there is a map with at most $c_3G^2$ crossing
  singularities.
\end{abstract}

\maketitle

\tableofcontents

\section{Introduction}

Among the different ways to combinatorially represent 3-manifolds, two
of the most popular are triangulations and surgery on a link.  A
triangulation is very natural way to represent 3-manifolds, and any
other representation of a 3-manifold is easy to turn into a
triangulation.  On the other hand, although some 3-manifold invariants
may be computed directly from a triangulation (e.g., the Turaev-Viro
invariants), not all can be, and it is difficult to visualise the
combinatorial structure of a triangulation.

A more typical way to present a 3-manifold is via Dehn surgery on a
link.  In practice, there are simple descriptions of small manifolds
via surgery, and this is generally the preferred way of representing
manifolds. There are many more invariants that may be computed
directly from
a surgery diagram, like the Witten-Reshetikhin-Turaev invariants.  It
is easy to turn a surgery diagram into a triangulation of the
manifold~\cite{Weeks05:ComputationHyperbolic}.
But for the other direction, passing from triangulations to surgery
diagrams, there seems to be little known.  In particular, it is
an open question whether a surgery diagram must (asymptotically) be more
complicated than a triangulation.  For a more general setting of this
problem, consider that if all the surgery coefficients are integers,
a surgery diagram naturally gives a 4-manifold bounded by the
3-manifold.  This leads us to ask the central question of the paper:

\begin{question}\label{quest:bound}
  How efficiently do 3-manifolds bound 4-manifolds?
\end{question}

To make this question more precise, let us make some definitions.

\begin{definition}
  \label{def:delta-triang}
  A \emph{$\Delta$-complex} is the quotient of a disjoint union of
  simplices by identifications of their faces.
  (See~\cite[Section~2.1]{Hatcher02:AlgebraicTopology} for a complete
  definition.)  A $\Delta$-triangulation is a $\Delta$-complex whose
  underlying topological space is a manifold.
\end{definition}

\begin{definition}
  The \emph{complexity} of a piecewise-linear oriented
  $n$-manifold~$M^n$
  is the minimal number of $n$-simplices in a
  $\Delta$-triangulation of~$M$.
  \begin{equation}
    \label{eq:triang-def}
    C(M^n) = \min_{\text{Triang.\ $\Delta$ of $M$}}\text{\# of $n$-simplices
      in $\Delta$}
  \end{equation}
\end{definition}

\begin{rem}
  Since the second barycentric subdivision of a $\Delta$-triangulation
  is an ordinary simplicial triangulation, $C(M)$ would only change by
  at most a constant factor if we insisted that the triangulation be
  simplicial.
\end{rem}

\begin{definition}
  The \emph{3-dimensional boundary complexity function}~$G_3(k)$ is
  the minimal complexity such that every 3-manifold of complexity at
  most~$k$ is bounded by a 4-manifold of complexity at most~$G_3(k)$.
\end{definition}

We can think of $G_3(k)$ as a kind of topological isoperimetric
inequality. We can now give a concrete version of our original
Question~\ref{quest:bound}:

\begin{question}\label{quest:simplicial-bound}
  What is the asymptotic growth rate of $G_3$?
\end{question}

The first main result of this paper is that $G_3(k) = O(k^2)$.  More
precisely, we have
\begin{citingthm}[Theorem \ref{teo:bound}]
  If a 3-manifold $M$ has a $\Delta$-triangulation with
  $t$~tetrahedra, then
  there exists a 4-manifold $W$ such that $\partial W=M$ and $W$~has
  a $\Delta$-triangulation
  with $O(t^2)$ simplices. Moreover, $W$ has ``bounded geometry''. That
  is, there exists an integer $c$ (not depending on $M$ and $W$) such
  that each vertex of the triangulation of $W$ is contained in less
  than $c$ simplices.
\end{citingthm} 

The fact that $W$~has bounded geometry makes the resulting
representation nicer; in particular, to check whether the topological
space resulting from a triangulation is a manifold, you need to decide
whether the link of each simplex is a sphere.  This is easy for
dimension $n\leq3$, in NP for
$n=4$~\cite{Schleimer04:SphereRecognition}, unknown for
$n=5$, and undecidable for dimension $n >
5$~\cite{Markov58:InsolubilityHomeomorphy,Markov58:UnsolvabilityCertain,LC03:UnrecognizabilitySphere}.
In all cases, such complexity issues do not arise if the triangulation
has bounded geometry.

Note that there is an evident linear lower bound for $G_3(k)$,  since a
triangulation for a 4-manifold also gives a triangulation of its
boundary.

We also prove a number of other related bounds which do not directly
refer to 4-manifolds.  For instance, we have
the following bound in terms of surgery:

\begin{citingthm}[Theorem \ref{teo:hyperbolic-surgery}]
  A finite-volume hyperbolic 3-manifold with volume~$V$
  has a rational surgery diagram with $O(V^2)$
  crossings.
\end{citingthm}

Note that there may be an infinite number of 3-manifolds with volume
less than the bound~$V$, and likewise an infinite number of surgeries
on a given link diagram; but in both cases the manifolds come in
families with some structure.  In this case as well there is a linear
lower bound: there are at least
$V/v_{\oct}$ crossings in any surgery diagram for~$M$, where
$v_{\oct} \approx 3.66$ is the volume of a regular ideal hyperbolic octahedron.
A somewhat weaker lower bound was proved by
Lackenby~\cite{Lackenby04:VolumeHypAlternating}; the bound using
$v_{\oct}$ comes from a decomposition into ideal octahedra, one per
crossing~\cite{Thurston99:HypVolume,Murakami01:KashaevInvariant}.

For a clean statement about general 3-manifolds, we use the
crucial notion of \emph{shadows}, which we recall in
Section~\ref{sec:shadows}.  For now, we just need to know that shadows
are certain kinds of decorated 2-complexes
which can be used to represent both a 4-manifold
and a 3-manifold (on the boundary of the 4-manifold), and that a
coarse notion of the
complexity of a shadow is the number of its \emph{vertices}.  There
are an infinite number of 3-manifolds with shadows with a given number of
vertices, but as with hyperbolic volume and surgeries on a given link,
they come in families that can be understood.
The \emph{shadow complexity} $\scomp(M)$ of a 3-manifold~$M$ is the
minimal number
of vertices in any shadow for~$M$.

The following theorem says that the shadow complexity gives a coarse
estimate of the hyperbolic volume.

\begin{citingthm}[Theorems \ref{teo:gromov-norm-lower} and \ref{teo:gromov-norm-upper}]
  There is a constant $C > 0$ so that every geometric 3-manifold~$M$,
  with boundary empty or a union of tori, satisfies
  \[
  \frac{v_{\tet}}{2v_{\oct}}\|M\| \le \scomp(M) \le C\|M\|^2.
  \]
  The lower bound on $\scomp(M)$ holds for all 3-manifolds.
\end{citingthm}

Here $v_{\tet}\approx 1.01$ is the volume of a regular ideal hyperbolic
tetrahedron and $v_{\oct}$ is as above.  A \emph{geometric manifold}
is one that satisfies the
Geometrization Conjecture~\cite{Thurston82:GeometryTopology}: it can
be cut along spheres and tori into pieces admitting a geometric
structure.  $\|M\|$~is the Gromov norm of~$M$, which is defined for any
3-manifold, and for a geometric 3-manifold is $1/v_{\tet}$ times
the sum of the volumes of the hyperbolic pieces.

Note that there is no constant term in these theorems.  The
manifolds with shadows with no vertices are the \emph{graph
  manifolds}, the geometric manifolds with no hyperbolic pieces (see
Proposition~\ref{prop:complexity-zero}).

Our techniques are based on maps from 3-manifolds to
surfaces, so we can also phrase the bounds in terms of the singularities
of such maps.  A \emph{crossing singularity} is a singularity of the
type we will consider
in Section~\ref{sec:extend-shadow}: a point in~$\RR^2$ with two
indefinite fold points in its inverse image.
For more background on the classification of
the stable singularities of a map from a 3-manifold to a 2-manifold, see
Levine~\cite{Levine85:ClassifyingImmersions,Levine88:Stable}.

\begin{citingthm}[Theorems \ref{teo:singularities-lower} and \ref{teo:singularities-upper}]
  A 3-manifold~$M$ has at least
  $\|M\|/10$ crossing singularities in any smooth, stable map $\pi:M\to\RR^2$.
  There is a universal constant~$C$ so
  that if $M$ is hyperbolic 
  then $M$ has a map to $\RR^2$ with $C\|M\|^2$ crossing singularities.
\end{citingthm}

One related theorem was previously known:
Saeki~\cite{Saeki96:SimpleStableMaps} showed that the manifolds with
maps to a surface with no crossing singularities are the graph
manifolds.

We also prove bounds for the complexity for 3-manifolds to bound
special types of 4-manifolds.

\begin{citingthm}[Theorem \ref{teo:bound-simply-conn}]
  A 3-manifold with a triangulation with $k$ tetrahedra is the
  boundary of a simply-connected 4-manifold with $O(k^2)$ 4-simplices.
\end{citingthm}

\begin{citingthm}[Theorem \ref{teo:bound-spin}]
  A 3-manifold with a triangulation with $k$ tetrahedra is the
  boundary of a spin 4-manifold with $O(k^4)$ 4-simplices.
\end{citingthm}

All the constants in these theorems can be made
explicit, but since in general they are quite bad, we have not
given them explicitly.  The exception is
Theorems~\ref{teo:gromov-norm-lower}
and~\ref{teo:singularities-lower}, which are
the best possible.

As one application of the results above, let us mention computing
invariants of 3-manifolds.  There are a number of 3-manifold
invariants that are most easily computed from a 4-manifold with
boundary.  (Often this is done via surgery diagrams, so the 4-manifold
is simply-connected, but there are usually more general constructions
as well.)  For instance, the Witten-Reshitikhin-Turaev (WRT) quantum
invariants are of this
form~\cite{RT91:Invariants3Mflds,Turaev94:QuantumInvariants,Roberts95:SkeinTheoryTV}
as is the Casson invariant~\cite{Lescop96:GlobalSurgeryCasson}.%
\footnote{The original definition of the Casson invariant is
  3-dimensional, but to compute it in practice the surgery formula is
  much easier.}

As one concrete example, Kirby and Melvin
explained~\cite{KM90:EvaluationsWRT,KM04:LocalSurgeryFormulas} how to
compute the WRT invariant at a 4th root of unity as a sum over spin
structures, which can be done concretely given a surgery diagram.
Although they show that the exact evaluation is NP-hard, our results
imply that the sum can be approximated (up to some error) in
polynomial time using random sampling over spin structures: for any
given spin structure, we can, in polynomial time, find a 4-manifold
which spin-bounds the given 3-manifold and therefore compute the
summand at this spin structure.  This contrasts with a result of
Kitaev and Bravyi, who showed that computing (up to the same error)
the partition function of the corresponding 2-dimensional TQFT is
BQP-complete, as soon as we allow evaluation of observables on closed
curves~\cite{BK00:QuantumInvariants}.

\medskip\noindent\textbf{Acknowledgments} We would like to warmly
thank Riccardo
Benedetti, Simon King, Robion Kirby, William Thurston, Vladimir
Turaev, and an anonymous referee for their encouraging comments and
suggestions.

\subsection{Plan of the paper}
\label{sec:plan-paper}

In the remainder of the Introduction, we survey some related work,
first on different kinds of topological isoperimetric functions, and
second on other work considering our main tool, stable maps from a
3-manifold to~$\RR^2$.  In Section~\ref{sec:4-mfld-stable-map}, we
sketch our construction in the smooth setting and introduce the
crucial tool of the Stein factorization, which shows how 2-complexes
naturally arise.  This section is not logically necessary for the rest
of the paper, although it does provide helpful motivation and a guide
to the proof.  This
2-complexes that arise are studied more abstractly in
Section~\ref{sec:shadows}, where we review the definition of shadow
surface and prove a number of properties of them; here we also use the
Gromov norm to prove all the lower bounds in the theorems above.  In
Section~\ref{sec:shadows-triangulations} we give our main tool, a
construction of a shadow from a triangulated 3-manifold with a map to
the plane, together with a bound on the complexity of the resulting
shadow.  Section~\ref{sec:shadows-triangulations} is independent from
Section~\ref{sec:shadows} except for the definition of shadows, and
only uses Section~\ref{sec:4-mfld-stable-map} as motivation, so the
impatient reader can skip there.  In Section~\ref{sec:upper-bounds} we
use this construction to prove the upper bounds of our main theorems
(except for the spin bound case, Theorem~\ref{teo:bound-spin})
and see precisely how shadow complexity relates to geometric notions
on the complexity of the manifold.  Finally, in
Section~\ref{sec:spin-boundaries} we show how to modify an arbitrary
shadow to get a 4-manifold that spin-bounds a specified spin structure
on a 3-manifold, while controlling the complexity.

\subsection{Related questions}
\label{sec:related-problems}

Although the question we consider does not seem to have been
previously addressed, there has been related work.  Perhaps the
closest is the work on distance in the pants
complex and hyperbolic volumes.  The pants complex is closely related to
shadows; in particular, a sequence of moves of length~$n$ in the pants
complex can be turned into a shadow with $n$ vertices for a 3-manifold
which has two boundary components, so that the natural pants
decomposition of the boundary components
corresponds to the start and end of the sequence of moves.

\begin{teo}[Brock~\cite{Brock03:WPVolumes,Brock03:WPMapping}]
  Given a surface $S$ of genus $g \geq 2$, there are constants $C_1,C_2$ so
  that for every pseudo-Anosov map $\psi:S\to S$, we have
  \[
  C_1\|\psi\|_{\mathrm Pants} \leq \vol(T_\psi) \leq C_2\|\psi\|_{\mathrm Pants}
  \]
  where $T_\psi$ is the mapping torus of $\psi$ and $\|\psi\|_{\mathrm Pants}$ is
  the translation distance in the pants complex.
\end{teo}

By the relation between moves in the pants complex and shadows
mentioned above, this
shows that for 3-manifolds that fiber over the circle with fiber a
surface of fixed genus, shadow complexity is bounded above and below
by a linear function of the hyperbolic volume.  However, the constant
depends on the genus in an uncontrolled way.  Our result gives a
quadratic bound, but with an explicit constant not depending on the
genus.  Brock's construction also produces shadows (and 4-manifolds)
of a particular type.

More recently, Brock and Souto have
announced~\cite{BS05:HeegardSplittings} that there is a similar
bound for manifolds with a Heegaard splitting with a fixed genus.  In
our language, their result says that a hyperbolic manifold with a
strongly irreducible Heegaard splitting of genus~$g$ has a shadow
diagram where the number of vertices is bounded by a linear function
of the volume, with a constant of proportionality depending only on
the genus.  (The result is probably true without the assumption that
the Heegaard splitting is strongly irreducible, but the statement
becomes more delicate in the language of the pants complex and we have
not checked the details.)  Their method of proof does not produce any
explicit constants.

There has also been work on the question of polygonal curves in
$\RR^3$ bounding surfaces.

\begin{definition}
  The \emph{surface isoperimetric function}~$G_{\text{\rm{surf}}}(k)$
  is the minimal number such that every closed polygonal
  curve~$\gamma$ in~$\RR^3$ with at most~$k$ segments bounds an
  oriented polygonal surface~$\Sigma$ with at
  most~$G_{\text{\rm{surf}}}(k)$ triangles.
\end{definition}

\begin{theorem}[Hass-Lagarias~\cite{HL03:MinimalNumber}]\label{teo:surf-bdry}
  $\frac{1}{2} k^2 \leq G_{\text{\rm surf}}(k) \leq 7k^2.$
\end{theorem}

This result contrasts sharply with the situation when we ask for the
spanning surface~$\Sigma$ to be a disk.

\begin{definition}
  The \emph{disk isoperimetric function}~$G_{\text{\textrm{disk}}}(k)$ is
  the minimal number such that every closed polygonal
  curve~$\gamma$ in~$\RR^3$ with at most~$k$ segments bounds an
  oriented polygonal disk~$D$ with at
  most~$G_{\text{\textrm{surf}}}(k)$ triangles.
\end{definition}

\begin{theorem}[Hass-Snoeyink-Thurston~\cite{HST03:SizeSpanning}]
  $G_{\text{\rm{disk}}}(k) = e^{\Omega(k)}$.  That is, there is a
  constant~$C$ so that, for sufficiently large~$k$,
  $G_{\text{\rm{disk}}}(k) \geq e^{Ck}$.
\end{theorem}

\begin{theorem}[Hass-Lagarias-Thurston~\cite{HLT03:AreaInequalities}]
  $G_{\text{\textrm{disk}}}(k) = e^{O(k^2)}$.  That is, there is a
  constant~$C$ so that, for sufficiently large~$k$,
  $G_{\text{\rm{disk}}}(k) \leq e^{Ck^2}$.
\end{theorem}

Although there is a large gap between these upper and lower bounds,
both bounds are substantially larger than the bounds in
Theorem~\ref{teo:surf-bdry}, which was about arbitrary oriented
surfaces.

There is an analogous question on the growth of $G_{\text{\rm{disk}}}$
for 3-manifolds rather than curves: asking for 4-balls bounding
a 3-sphere with a
given triangulation on the boundary.  As stated, this is not an
interesting question, since we can construct such a
triangulation by taking the triangulated 3-ball and coning it to a
point.  This is related to the somewhat unsatisfactory nature of the
4-manifold complexity (mentioned earlier).  A more interesting
question might involve 4-manifold triangulations where the vertices
have bounded geometry.  For a somewhat different question, there are
known upper bounds:

\begin{definition}
  The \emph{Pachner isoperimetric function}
  $G_{\text{\textrm{Pachner}}}(k)$ is the maximum over all
  triangulations~$T$ of the 3-sphere with $\leq\,k$ simplices of the
  minimum number of Pachner moves required to relate~$T$ to the
  standard triangulation, the boundary of a 4-simplex.
\end{definition}

\begin{theorem}[King~\cite{King01:Polytopality}, Mijatovi\'c~\cite{Mijatovic03:Simplifying}]
  $G_{\text{\rm{Pachner}}}(k) = e^{O(k^2)}$.
\end{theorem}

Note that a sequence of Pachner moves as in the definition gives you,
in particular, a triangulation of the 4-ball, although you only get
very special triangulations of the 4-ball in this way.

As in the case of polygonal surfaces and disks, this upper bound is
much larger than the polynomial bound we obtain.
King~\cite{King01:Polytopality} also constructs triangulations
of~$S^3$ which seem likely to require a large number of Pachner moves
to simplify.

\subsection{Previous work}
\label{sec:previous-work}

The central construction
in our proof, a generic smooth map from a 3-manifold to $\RR^2$, has
been considered by several previous authors, sometimes with little
contact with each other.  These maps were probably first considered by
Burlet and de Rham~\cite{BdR74:CertainsApplications}, who showed that
the 3-manifolds admitting a map with only definite fold singularities
are connected sums of $S^1 \times S^2$ (including $S^3$).  They also
introduced the Stein factorization.
Levine~\cite{Levine88:Stable} clarified the
structure of the singularities and studied, for instance, related
immersions of the 3-manifold into $\RR^4$.  Burlet and de Rham's
result was extended by Saeki~\cite{Saeki96:SimpleStableMaps}, who
showed that the 3-manifolds admitting a map without codimension 2
singularities (i.e., only definite or indefinite folds) are the graph
manifolds.

Rubinstein and Scharlemann~\cite{RS96:ComparingHeegardSplittings}
constructed a map from the complement of a graph in a 3-manifold
to~$\RR^2$ from a pair of Heegaard splittings and used this to bound
the number of stabilizations required to turn one splitting into the
other.  Much of the analysis is similar to ours.

Also independently, Hatcher and
Thurston~\cite{HT80:PresentationMappingClassGroup} considered Morse
functions on a orientable surface to show that its mapping
class group is finitely presented.  To get a set of generators, they
considered one parameter deformations of the Morse function. Note that
a one parameter family of maps from~$\Sigma$ to~$\RR$ is a map from the
3-manifold $\Sigma\times[0,1]$ to $\RR^2$.

In a slightly different context, Hatcher's proof of the Smale
conjecture~\cite{Hatcher83:ProofSmale}, that the space of smooth
2-spheres in $\RR^3$ is contractible, uses the Stein factorization of
a map from $S^2$ to $\RR^2$.  Hong, McCullough, and Rubinstein
recently combined this approach with the Rubinstein-Scharlemann
techniques in their proof of the Smale conjecture for lens
spaces~\cite{HMR04:SmaleConjectureLens}.

On the other side of the story,
Turaev~\cite{Turaev91:TopologyShadows,Turaev92:ShadowLinks,Turaev94:QuantumInvariants}
introduced shadow surfaces as the most natural objects on which the
Reshetikhin-Turaev quantum invariants are defined.  He observed that
you could construct both a 3-manifold and a 4-manifold with boundary
the 3-manifold from a shadow surface.

Thus several authors have been considering nearly the same objects
(shadow surfaces on the one hand and the Stein factorization of a map
from $M^3$ to $\RR^2$ on the other hand) for several years.  The gleams
are key topological data
from the shadow surface point of view, since they let you reconstruct
the 3-manifold, but they were not explicitly described by the authors
writing on Stein factorizations, although it is implicitly present.

\section{4-manifolds from stable maps}
\label{sec:4-mfld-stable-map}

In order to prove that 3-manifolds efficiently bound 4\hyp manifolds,
we start by sketching a proof that 3-manifolds do bound 4-manifolds.
In Section~\ref{sec:shadows-triangulations} we will analyze a version of the
proof in the PL setting and give a bound on the complexity of the
resulting 4-manifold.

Consider an oriented, smooth,
closed 3-manifold~$M^3$ and a generic smooth map~$f$ from $M$ to
$\RR^2$.
At a regular value $x \in \RR^2$, the inverse image $f^{-1}(x)$ consists
of an oriented union of circles.  To construct a
4-manifold, we glue a disk to each of these circles away from
critical values and then extend across the singularities in
codimension 1 and codimension 2.

\subsection{Pants decompositions from Morse functions}
\label{sec:pants-decomp}

\begin{figure}[htbp]
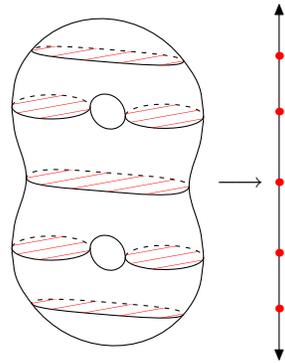

  \[ \mfigb{draws/surface.62} \longrightarrow \mfigb{draws/surface.71} \]
  \caption{Proving that every surface bounds a 3-manifold: A generic
    map of a surface to $\RR$, with regular values marked and disks
    glued into the inverse image of the regular values.}
  \label{fig:cut-surface}
\end{figure}

To get some idea of what the singularities look like, we first do the
analysis of extending across
singularities one dimension down: let's
prove that every oriented 2-manifold~$\Sigma^2$ bounds a 3-manifold.
Consider a generic smooth map~$f$ from~$\Sigma$ to~$\RR$, that is, a Morse
function.  The inverse image of a regular
value is again a union of circles.  Glue in disks to each of these circles
as in Figure~\ref{fig:cut-surface}.  More properly, take $\Sigma\times[0,1]$,
pick a regular
value in each component of $\RR$ minus the singular set, and attach
2-handles along the circles appearing in the
inverse image of the chosen regular values.  The result is a
3-manifold with one boundary component which is~$\Sigma$ and other boundary
components corresponding to the singular values of~$f$.

The singular values of a Morse function, locally in the domain~$\Sigma$,
are well-known: they are critical points with a quadratic form which
is definite (index 0 or 2, minima or maxima) or indefinite (index 1, saddle
points).  Since our construction works with the entire
inverse image of a regular value, we need to understand the
singularities locally in the range~$\RR$; that is, we need to know the
connected components of inverse images of a critical value.  This is
easy for the definite singularities.

Let $p_0\in\Sigma$ be a saddle point, and let $x_0\in\RR$ be its
image.  Near $p_0$, $f^{-1}(x_0)$ is a cross.  For $x$ above and below
$x_0$, $f^{-1}(x)$ is locally the cross is smoothed out in the two
possible ways. Note that
the orientations of~$\Sigma$ and~$\RR$ induce an orientation of~$f^{-1}(x)$
for all~$x\in\RR$ except at critical points in $\Sigma$, so both of these
smoothings must be oriented, so the arms of the
cross must be oriented alternating in and out.  The connected
component of $f^{-1}(x_0)$ containing $p_0$ must join the arms of the
cross in an orientation-preserving way and is therefore a figure 8 graph~$\mfigb{draws/surface.100}$.
This implies that for a small interval~$I$ containing~$x_0$,
$f^{-1}(I)$ is a pair of pants.

\begin{figure}[htbp]
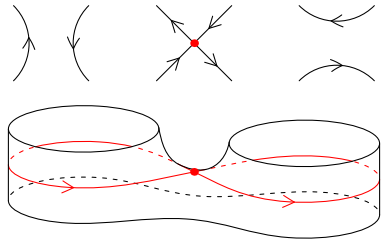

  \[ \mfigb{draws/surface.5}\qquad
     \mfigb{draws/surface.0}\qquad
     \mfigb{draws/surface.6} \]
  \[ \mfigb{draws/surface.15} \]
  \caption{Analysing the saddle singularity}
  \label{fig:saddle-sing}
\end{figure}

To finish constructing the 3-manifold, recall that in the previous step
we glued in disks at all the regular values.  Near this pair of pants,
this means that we have closed off each hole in the pair of pants and
the boundary component we are trying to fill in is just a sphere,
which we can fill in with a ball.

An easier analysis shows that the surface we need to fill in the
other cases of a maximum or minimum is again a sphere.

Notice where the proof breaks down if we do not assume that~$\Sigma$
is oriented: there is then another possibility for the inverse image
of the critical value, with opposite arms of the cross attached to each
other: $\mfigb{draws/surface.101}$.  In this case the
surface we are left to fill in turns out to be $\mrp^2$, which does
not bound a 3-manifold.

In a similar way we can analyse the possible stable singularities of a
smooth map
from a 3\hyp manifold~$M$ to $\mr^2$.  We glue in a disk (a 2-handle) to
each circle in the inverse image of a regular point, extend across
codimension 1 singularities by attaching 3-handles (the singularities
look just like the singularities we analysed for the case of a
surface, crossed with $\RR$), and then consider the codimension~2
singularities.  In Section~\ref{sec:extend-shadow}, we will analyze
the codimension~2 singularities (in the PL category) and see that the
remaining boundary from each codimension~2 singularity is~$S^3$, which
can be filled in by attaching a 4-handle.

\subsection{Stein factorization and shadow surfaces}
\label{sec:stein-surface}

\begin{figure}[htbp]
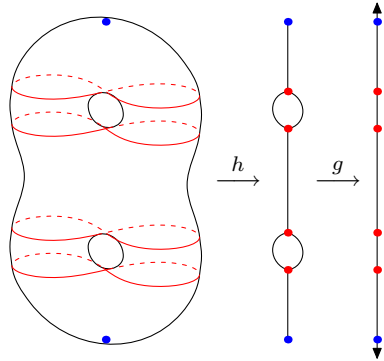

  \[ \mfigb{draws/surface.63} \overset{h}{\longrightarrow} \mfigb{draws/surface.75} \overset{g}{\longrightarrow}
  \mfigb{draws/surface.72}
  \]
  \caption{The Stein factorization $f = g\circ h$ of the map in Figure~\ref{fig:cut-surface}.}
  \label{fig:stein-graph}
\end{figure}

For a more global view, we can consider the \emph{Stein factorization}
$f = g \circ h$ of these maps.  The Stein factorization of a map~$f$ with
compact fibers decomposes it as the composition of a map~$h$ with
connected fibers and a map~$g$ which is finite-to-one.  That is, $h$
is the quotient onto the space of connected components of the fibers
of~$f$.  See Figure~\ref{fig:stein-graph} for an example.

For a stable map from an oriented surface~$\Sigma$ to $\RR$, the Stein
factorization is generically a 1-manifold, with singularities from the
critical points.  Concretely, it is a graph with vertices which have
valence 1 (at definite singularities) or valence 3 (at indefinite
singularities).  The
surface is a circle bundle over this Stein graph~$\Gamma$ at generic points.
Likewise, the 3\hyp manifold we constructed (with boundary~$\Sigma$) is a
disk bundle over~$\Gamma$ at generic points.  In fact, the 3-manifold
collapses onto~$\Gamma$.
If we collapse all the valence~1 ends, we may think of~$\Gamma$ as
representing a pair-of-pants decomposition of~$\Sigma$; each circle in
the pair-of-pants decomposition bounds a disk in the 3\hyp manifold.

\begin{figure}[htbp]
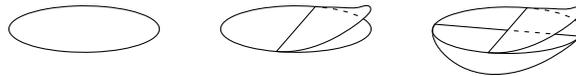

  \[\mfig{draws/surface.90}\qquad\mfig{draws/surface.91}\qquad\raisebox{-3.9bp}{\hbox{$\mfig{draws/surface.92}$}}\]
  \caption{Some local models for the Stein factorization of a map from
    a 3-manifold to $\RR^2$ in codimension 0, 1, and 2.  In each
    picture, the map to the plane is the vertical projection.}
  \label{fig:stein-surface}
\end{figure}

For a stable map from a 3-manifold to $\RR^2$, on the other hand, the
Stein factorization is generically a surface.  The
codimension~1 singularities of the Stein surface are products of the
lower-dimensional singularities with an interval, and have one
or three sheets meeting at an edge at what we will call definite or
indefinite folds, respectively.  In codimension two there are a
few different configurations of how the surface can meet, the most
interesting of which is shown in Figure~\ref{fig:stein-surface}.

The 3-manifold is a circle bundle over the Stein surface at generic
points and the 4-manifold is generically a disk bundle.  As in the
previous case, it turns out that the 4-manifold collapses onto the
Stein surface.  The resulting surface is very close to a shadow surface.

Unlike in the lower dimensional case, the surface does not determine
the 4-manifold (or the 3-manifold), even after you fix a standard
local model of how the surface sits inside the 4-manifold.  The
additional data you need are the \emph{gleams}, numbers associated to
the 2-dimensional regions of the surface; see
Definition~\ref{def:gleams}.

\section{Shadow surfaces}
\label{sec:shadows}

We will now define shadow surfaces and shadows of
3-manifolds and give a few examples.   In the
Section~\ref{sec:relativeshadows} we
will extend the definitions to 3-manifolds with boundary and an
embedded, framed graph.  For a
more detailed though introductory account of shadows of 3-
and 4-manifolds, see~\cite{Costantino04:ShortIntroduction}.  Note that
these are
slightly different from the Stein surfaces mentioned in
Section~\ref{sec:stein-surface}. We prefer shadow surfaces as the
fundamental object since they are a little more symmetric and
regular than Stein surfaces. From now on every manifold will be PL
compact and oriented unless explicitly stated and every polyhedron
will be finite; we also recall that, in dimension 3 and 4,
each PL manifold has a unique smooth structure and vice versa.

\subsection{Shadows of 3-manifolds}
\label{sec:shadows-3-manifolds}

For simplicity, we will first define shadows in the case when there is
no boundary, appropriate for 3\hyp manifolds without boundary or other
decorations; in the next section we will extend this.

\begin{figure}
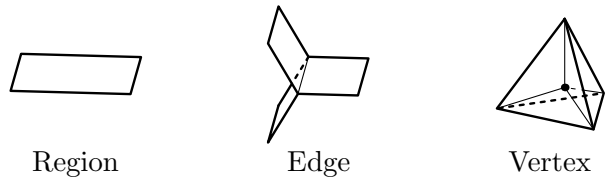

\[
\begin{array}{c@{\qquad\qquad}c@{\qquad\qquad}c}
\mfigb{draws/surface.110}&\mfigb{draws/surface.111}&\mfigb{draws/surface.112}\\
\textrm{Region}&\textrm{Edge}&\textrm{Vertex}
\end{array}
\]
  \caption{For each point of a simple polyhedron embedded in a 4-manifold there is a local chart
    $(U,\phi)$ in which the polyhedron is flat in the sense that is
    embedded in a three dimensional plane and in this plane it appears
    as in one of the three models shown in this figure.  }
  \label{fig:singularityinspine}
\end{figure}

\begin{defi} A simple polyhedron $P$ is a compact topological space
  where every point has a neighborhood homeomorphic to an open set in
  one of the local models depicted in
  Figure~\ref{fig:singularityinspine}.
  The set of points without a local model of the leftmost type form a
  4-valent graph, called the
  \emph{singular set} of the polyhedron and denoted~$\Sing(P)$. 
  The vertices of $\Sing(P)$ are called \emph{vertices} of~$P$.  The
  connected components of $P\setminus\Sing(P)$ are called the
  \emph{regions} of $P$.  The set of points of $P$ whose local models
  correspond to the boundaries of the blocks shown in the figure is
  called the {\it boundary} of~$P$ and is denoted~$\partial P$; $P$
  is said to be \emph{closed} if it has empty boundary. A region is
  {\it internal} if its closure does not touch $\partial P$.  A simple
  polyhedron is \emph{standard} if every region of $P$ is a disk, and
  hence $\Sing(P)$ has no circle components.
\end{defi}

\begin{defi}
  Let $W$ be a PL, compact and oriented 4-manifold.  $P \subset W$ is a
  \emph{shadow} for $W$ if~$P$ is a closed simple sub-polyhedron onto which
  $W$ collapses and~$P$ is \emph{locally flat} in $W$, that is for each point
  $p\in P$ there exists a local chart $(U,\phi)$ of $W$ around $p$ such
  that $\phi(P\cap U)$ is contained in $\RR^3 \subset \RR^4$.
\end{defi}

It follows from this definition that in the 3-dimensional slice, the
pair $(\mathbb{R}^3\cap \phi (U),\allowbreak\mathbb{R}^3\cap \phi(U\cap P))$ is
PL-homeomorphic to one of the models depicted in
Figure~\ref{fig:singularityinspine}.

For the sake of simplicity, from now on we will skip the PL prefix
and all the homeomorphisms will be PL unless explicitly stated.
Not every 4-manifold admits a shadow: a necessary and sufficient
condition for $W$ to admit one is that it has an handle
decomposition containing no handles of index greater
than~2~\cite{Turaev91:TopologyShadows,Costantino04:ShadowsBranchedShadows}.
This imposes restrictions on the
topology of~$W$.  For instance, its boundary has to be a non-empty
connected 3-manifold. 

\begin{defi}
A shadow of an oriented, closed 3-manifold $M$ is a shadow $P$ of a closed 4-manifold~$W$
with $M = \partial W$.
\end{defi}
\begin{teo}[Turaev~\cite{Turaev91:TopologyShadows}]
Any closed, oriented, connected 3-manifold has a shadow.
\end{teo}

\begin{example}\label{examp:spheres}
  The simple polyhedron $P=S^2$ is a shadow of $S^2\times D^2$ and
  hence of the 3-manifold $S^2\times S^1$. In this case, $P$ is a
  surface whose self-intersection number in the ambient 4-manifold is
  zero. Consider now the disk bundle over $S^2$ with Euler number
  equal to 1, homeomorphic to a punctured $\mathbb{CP}^2$. The
  0-section of the bundle is a shadow of the 4-manifold homeomorphic
  to $P$ and so $P$ is a shadow of $\mathbb{CP}^2-B^4$ and of its
  boundary: $S^3$.
\end{example}

The above example shows that the naked polyhedron by itself is not
sufficient to encode the topology of the 4-manifold collapsing on it.
Turaev described~\cite{Turaev91:TopologyShadows} how to equip a
polyhedron embedded in a 4-manifold with combinatorial data called
{\it gleams} which are sufficient to encode the topology of the
regular neighborhood of the polyhedron in the manifold. A gleam is a
coloring of the regions of the polyhedron with values in
$\frac{1}{2}\ZZ$, with value modulo~1 given by a $\ZZ_2$-gleam which
depends only on the polyhedron.

In the simplest case, if $P$ is a shadow of $M$ and $P$ is
homeomorphic to an orientable
surface, then $W$ is homeomorphic to an oriented disk bundle over the
surface and the gleam of $P$ is the Euler number of the normal bundle
of $P$ in $W$.

We summarize in the following proposition the basic construction
of the $\mathbb{Z}_2$-gleam and of the gleam of a simple
polyhedron. A \emph{framing} for a graph~$G$ in a 3-manifold~$M$ is a
surface with boundary embedded in~$M$ and collapsing onto~$G$.
\begin{prop}\label{prop:embpoly}
  Let $P$ be a simple polyhedron. There exists a canonical
  $\mathbb{Z}_2$-coloring of the internal regions of $P$ called the
  \emph{$\mathbb{Z}_2$-gleam} of $P$. If $P$ is embedded in a
  4-manifold $W$ in a flat way, there is a canonical coloring of the
  internal regions of $P$ by integers or half-integers called
  \emph{gleams}, such that the gleam of a region of~$P$ is an integer
  if and only if its $\mathbb{Z}_2$-gleam is zero. Moreover, if
  $\partial P\subset \partial W$ is framed, then the gleam can also be
  defined on the non-internal regions of $P$.
\end{prop}
\begin{prf}
Let $D$ be an internal region of $P$ and let
$\overline{D}$ be the abstract compactification of the (open) surface
represented by $D$. The embedding of $D$ in $P$ extends to a map
$i:\overline{D}\to P$ which is injective on $\intr(\overline{D})$,
locally injective on $\partial \overline{D}$ and sends $\partial
\overline{D}$ into $\Sing(P)$. Using $i$ we can ``pull back'' a small
open regular neighborhood of $D$ in $P$ and construct a simple
polyhedron $U(D)$ collapsing on $\overline{D}$. Extend $i$ as a local
homeomorphism $i':U(D)\to P$ whose image is contained in a small
regular neighborhood of the closure of $D$ in $P$. In the particular
case when $i$ is an embedding of $\overline{D}$ in $P$, $U(D)$ is the
regular neighborhood of $D$ in $P$ and $i'$ is its embedding in
$P$. In general, $U(D)$ has the following structure: each boundary
component of $\overline{D}$ is glued to the core of a band (annulus or
M\"obius strip) and some small disks are glued along half of their
boundary on segments which are properly embedded in these bands and
cut transversally once their cores. We define the $\mathbb{Z}_2$-gleam
of $D$ in $P$ to be equal to the reduction modulo 2 of the number of
M\"obius strips used to construct $U(D)$. This coloring only depends
on the combinatorial structure of $P$.

Let us now suppose that $P$ is embedded in a 4-manifold $W$, and let
$D$, $\overline{D}$, $i:\overline{D}\to P$, $U(D)$ and $i'$ be defined
as above. Using~$i'$, we can ``pull back'' a neighborhood of
$\overline{D}$ in~$W$ to an oriented 4-ball $B^4$ collapsing on
$U(D)$. The regular neighborhood of a point $p_0\in \partial
\overline{D}\subset U(D)$ sits in a 3-dimensional slice~$B^3_0$
of~$B^4$ where it appears as in Figure~\ref{fig:divergingdirection}.
The direction along which the other regions touching $\partial
\overline{D}$ get separated gives a section of the bundle of
orthogonal directions to~$D$ in~$B^4$.
(If $p_0\in \partial P$, use the framing of~$\partial P$ in place of
the other regions.  By an \emph{orthogonal direction} we mean a line
in the normal bundle, not a ray.)  This section can be defined on
all $\partial \overline{D}$ and the obstruction to extend it to all
of~$\overline{D}$ is an element of $H^2(\overline{D},\partial
\overline{D};\pi_1(S^1))$. Since~$B^4$ is oriented, we can
canonically identify this element with an integer and define the gleam
of~$D$ to be half this number.
\begin{figure}
  \centerline{\includegraphics[width=6.4cm]{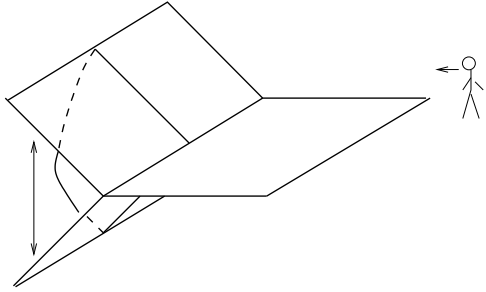}}
  \caption{The picture sketches the position of the polyhedron in a
    3-dimensional slice of the ambient 4-manifold.  The horizontal
    plane is our region of interest~$D$.  The direction
    indicated by the vertical double arrow is the one along which the
    two regions touching the horizontal one get separated.  }
  \label{fig:divergingdirection}
\end{figure} 
\end{prf}

We can also go the other direction, from gleams to 4-manifolds.
\begin{defi}\label{def:gleams}
  A gleam on a simple polyhedron $P$ is a coloring on all the regions
  of $P$ with values in $\frac{1}{2}\mathbb{Z}$ such that the color of
  each internal region is integer if and only if its
  $\mathbb{Z}_2$-gleam is zero.
\end{defi}
\begin{teo}[Reconstruction of 4-manifold~\cite{Turaev91:TopologyShadows}]\label{teo:reconstruction}
  Let $P$ be a polyhedron with gleams~$g$; there exists a canonical
  reconstruction associating to $P$ and $g$ a pair $(W_P,P)$ where
  $W_P$ is a PL, compact and oriented 4-manifold containing a properly
  embedded copy of $P$ with framed boundary, onto which it collapses
  and such that the gleam of $P$ in $W_P$ coincides with $g$. The pair
  $(W_P,P)$ can be explicitly reconstructed from the combinatorics of
  $P$ and from its gleam. Moreover, if $P$ is a polyhedron embedded in
  a PL and oriented manifold $W$, $\partial P$ is framed, and $g$ is
  the gleam induced on $P$ as explained in the
  Proposition~\ref{prop:embpoly}, then $W_P$ is homeomorphic to the
  regular neighborhood of $P$ in $W$.
\end{teo}

Hence, to study 4-manifolds with shadows (and their boundaries), one
can either use abstract polyhedra equipped with gleams or embedded
polyhedra.

From now on, each time we
speak of a shadow of a 3-manifold as a polyhedron we will be
implicitly taking a 4-dimensional thickening of this polyhedron
whose boundary is the given 3-manifold or, equivalently,
a choice of gleams on the regions of the polyhedron.
\begin{example}
  Let $M$ be a 3-manifold which collapses onto a simple polyhedron $P$
  whose regions are orientable surfaces; it is straightforward to
  check that the $\mathbb{Z}_2$-gleam of $P$ is everywhere zero.  Let
  us then equip $P$ with the gleam which is zero on all the regions;
  Turaev's thickening construction produces the 4-manifold $W=M\times
  [-1,1]$ and $P$ is a shadow of~$\partial W$, homeomorphic to the
  double of~$M$.
\end{example}

\subsection{Shadows of framed graphs in manifolds with boundary}\label{sec:relativeshadows}
In this subsection we extend the definition of shadows
to pairs $(M,G)$ where $M$~is a
3-manifold, possibly with boundary and $G \subset M$ is a (possibly empty)
framed graph with trivalent vertices and univalent ends.  To do this,
we allow the simple polyhedron to have boundary.
\begin{defi}\label{def:bndry-decorated}
  A \emph{boundary\hyp decorated} simple polyhedron~$P$ is a simple polyhedron
  where  $\partial P$
  is equipped with a cellularization whose 1-cells are colored with one
  of the following colors: $i$ (\emph{internal}), $e$ (\emph{external})
  and $f$ (\emph{false}). We can correspondingly distinguish three
  subgraphs of $\partial P$ intersecting only in $0$-cells and whose
  union is~$\partial P$: let us call them $\partial_i P$, $\partial_e
  P$, and $\partial_f P$.  A boundary\hyp decorated simple polyhedron
  is said to be \emph{proper} if $\partial_f(P) = \emptyset$.
\end{defi}

We can turn decorated polyhedra into shadows.  The intuition is that
$\partial_f(P)$ is ignored, $\partial_e(P)$ is drilled out to create
the boundary of the 3-manifold, and $\partial_i(P)$ gives a trivalent
graph.

\begin{defi}
  Let $P$ be a boundary\hyp decorated simple polyhedron, properly embedded in a
  4-manifold~$W$ which collapses onto~$P$ with a framing on
  $\partial_i(P)$.  Let $M$~be the complement of an open regular
  neighborhood of~$\partial_e P$ in~$\partial W$, and let $G$ be a
  framed graph whose core is $\partial_i(P)$.  Then we say that~$P$ is
  a \emph{shadow} of $(M,G)$ and, if~$\partial_f P=\emptyset$, we call
  it a \emph{proper shadow}.  As before, we can define a \emph{gleam}
  on each region of~$P$ that does not meet $\partial_f(P)
  \cup \partial_e(P)$.
\end{defi}
\begin{rem}
  Turaev's Reconstruction Theorem extends to the case of decorated
  polyhedra equipped with gleams on the regions not touching
  $\partial_e P\cup\partial_f P$.
\end{rem}
\begin{rem}
  \label{rem:genus-boundary}
  The genus of a boundary component of~$M$ equals the rank of $H_1$ of
  the corresponding component of $\partial_e(P)$, since the Euler
  characteristic of the handlebody filling the boundary component
  equals the Euler characteristic of the graph.  In particular,
  components of~$\partial_e(P)$ corresponding to sphere boundary
  components of~$M$ are contractible and so if
  $\partial_f(P)$ is empty, $M$ does not have any sphere boundary
  components that do not meet~$G$.
\end{rem}
\begin{teo}[Turaev~\cite{Turaev91:TopologyShadows}]
  Let $M$ be a oriented, connected 3-manifold and let~$G$ be a
  properly-embedded
  framed graph in~$M$ with vertices of valence~1 or~3.  If $M$ has no
  spherical boundary components that do not meet~$G$, the pair
  $(M,G)$ has a proper, simply-connected shadow.
\end{teo}

Let us now show how to construct a shadow of a pair $(M,G)$ given a
shadow~$P$ of $(M,\emptyset)$. Recall that $M$ is the boundary of a
4\hyp manifold collapsing through a projection~$\pi$ onto~$P$.
Up to small isotopies, we can suppose
that the restriction to~$G$ of~$\pi$ is
transverse to $\Sing(P)$ and to itself; that is, it does not contain
triple points or self tangencies and is injective on the vertices of
$T$.  Let us also suppose that it misses $\partial_f(P)$.
Then the mapping cylinder of the projection of $G$ in $P$ is
contained in the thickening~$W_P$ of~$P$ and $W_P$ collapses on it.
(Recall that the mapping cylinder is $P \cup G \times [0,1]$, with $G
\times\{0\}$ identified with $\pi(G) \subset P$.)
By Proposition~\ref{prop:embpoly} we can equip this polyhedron with
gleams. Coloring $G\times \{1\} $ with the color $i$ we get a shadow
of the pair $(M,G)$, coloring it with $e$ we get a shadow of $M-U(G)$
where $U(G)$ is a small open regular neighborhood of $G$ in $M$, and
coloring it with $f$ we get another shadow of $M$ (necessarily not
proper).

As a warm up, note that a flat disk~$D$ whose
boundary has color~$f$ is a shadow of the pair $(S^3,\emptyset)$.  The
open solid torus~$T_h=\pi^{-1}(\intr(D))$ can be imagined as the
regular neighborhood of the closure of the $z$-axis in $\mr^3$,
embedded inside~$S^3$ in the standard way.  The fibers of $\pi:S^3\to D$ run
parallel to the $z$-axis away from~$\infty$ and are unknotted. The
projection of the solid torus $T_v=S^3 \setminus T_h$ is $\partial D$.

With the setup above, we now apply the projection construction to the
case of a link~$L$ in~$S^3$.  Up to isotopy we can suppose that
$L\subset T_h$ and that its projection to $D$ is generic; so it is
sufficient to consider a standard diagram of $L$ in the unit disk in
$\mr^2$.  The mapping cylinder~$D_L$ of $\pi:L\to D$ is obtained
from~$D$ by gluing an annulus for each component of~$L$ and marking
the free boundary components of these annuli with the color~$i$.  We can further
collapse the region of~$D_L$ containing $\partial_f D$; this produces
a simple sub-polyhedron of~$D_L$, which we call~$P_L$.  By
construction~$\partial P_L=\partial_i P_L=L$.

\begin{figure}
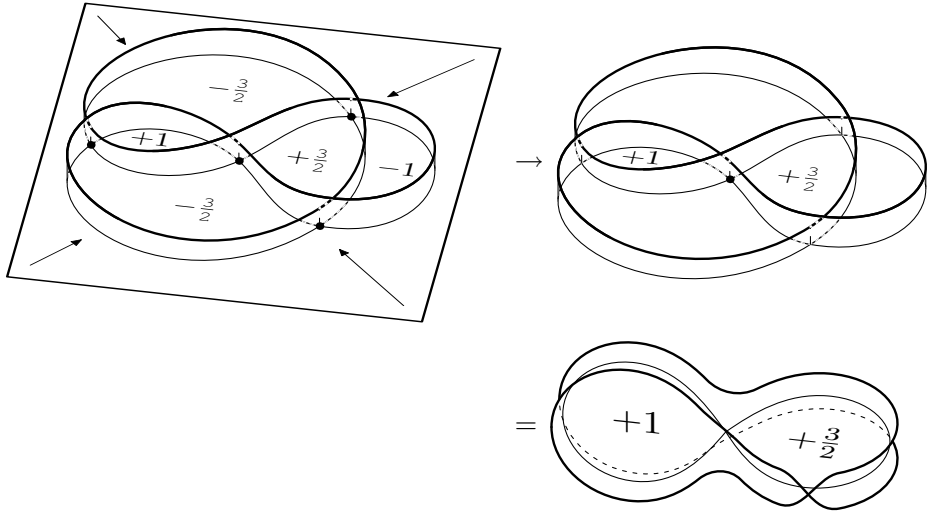

  \begin{align*}
  \mfigb{draws/shadow.1} &\to \mfigb{draws/shadow.2}\\
      &= \mfigb{draws/shadow.11}
%      &\cong \mfigb{draws/shadow.10}
  \end{align*}
  \caption{In the left part of the picture we sketch the construction
    described in Example~\ref{examp:figureeight}; the resulting shadow
    is drawn in the right part where, for the two internal regions we
    write their gleams. Note that, after the collapse of the
    polyhedron~$D_L$ along its free boundary component (as indicated
    by the arrows), the only vertex surviving is the central one. }
  \label{fig:figureeight}
\end{figure}

In general, $P_L$ has some vertices, each corresponding to a crossing
in the diagram of~$L$. However, some of the crossings in that
diagram of~$L$ do not generate vertices in $P_L$ because they
disappear when we pass from~$D_L$ to~$P_L$.
\begin{example}\label{examp:figureeight}
  Applying the construction to a figure eight knot in a standard
  position, one gets a shadow of its complement containing only one
  vertex: Three of the four crossings of the diagram are
  contained in the boundary of the region to be collapsed in~$P_L$.
  See Figure~\ref{fig:figureeight}.
\end{example} 

\begin{figure}
  \centering
  \includegraphics{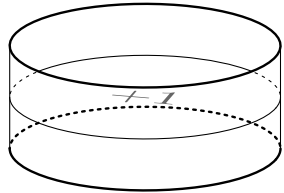}
  \caption{The shadow obtained for the Hopf link.}
  \label{fig:hopf}
\end{figure}

\begin{example}\label{examp:hopflink}
  Consider a standard Hopf link in $T_h$. The polyhedron $H$ one gets
  by applying the above procedure contains no vertices: the two
  crossings of the standard projection of the Hopf link in $\mr^2$
  touch the region of $D$ which is collapsed. The resulting polyhedron
  can be obtained by gluing a disk to the core of an annulus; its
  embedding in $B^4$ is such that $\partial H$ is the Hopf link in
  $\partial B^4$.  See Figure~\ref{fig:hopf}.
\end{example}

\begin{example}\label{examp:secondsingularity}
  Applying this construction to the graph~$\partial P$ in
  Figure~\ref{fig:embeddedsecondsingularity} (with the projection
  shown there) gives the shadow in Figure~\ref{fig:secondopoliedro}.
\end{example}

\subsection{Shadow complexity and its basic properties}
\label{sec:shadow-complexity-basic}
We will now define a notion of \emph{shadow complexity} and study how
it behaves under combining
manifolds, either via connect sum (gluing along spheres) or torus
connect sum (gluing along torus boundaries).

\begin{defi}
  For $M$ be an oriented 3-manifold (possibly with boundary) and
  $T\subset M$ a trivalent graph, the \emph{shadow
    complexity}~$\scomp(M,T)$ of the pair $(M,T)$ is the minimal
  number of vertices of a boundary-decorated shadow of $(M,T)$.
\end{defi}

\begin{rem}\label{rem:collapse-false-bdry}
If $M$ has no spherical boundary components,
it doesn't matter whether or not we allow the shadow to have false
edges in this definition: If we have a decorated shadow~$P_1$ for
$(M,T)$, it can be shown that the polyhedron~$P_2$ obtained by
iteratively collapsing all the regions of~$P_1$ containing a false
boundary edge is a complex obtained by gluing some graphs to a
(possibly disconnected) simple polyhedron.  This complex~$P_2$ can be
modified, without adding vertices, to give a shadow~$P_3$ for $(M,T)$
without false edges and no more vertices than~$P_1$; the modifications
are generally similar to those in Lemma~\ref{lem:sc-connectsum}, with
a few special constructions for cases where the complex is contractible
(so $M$ is $S^3$) or the graph has non-trivial loops, producing
$S^1\times S^2$ summands in the prime decomposition.  All of these
special cases are graph manifolds.  By
Proposition~\ref{prop:complexity-zero} they can be treated without
creating any vertices.
\end{rem}

Such a notion of complexity is similar to the usual notion
of complexity of 3-manifolds introduced by
S. Matveev~\cite{Matveev89:ComplexityThreeDimensional}:

\begin{defi}
  \label{def:complexity}
  The \emph{complexity}~$c(M)$ of a 3-manifold~$M$ is the minimal number of
  vertices in a simple polyhedron~$P$ contained in~$M$ which is a
  spine for $M$ or $M$ minus a ball.
\end{defi}

Both notions are based
on the least number of vertices of a simple polyhedron describing (in a
suitable sense) the given manifold.
Despite this similarity, shadow
complexity is not finite.  That is, the set of manifolds having complexity
less than or equal to any given integer is infinite. For
instance, the lens spaces $L(p,1)$ have a shadow surface which is
$S^2$ with gleam~$p$ and so they all have shadow complexity~0.

To reduce the set of attainable manifolds to a finite number and bound
the complexity of the
4-manifold, we also need to bound the gleams.

\begin{defi}\label{def:gleam-weight}
  The \emph{gleam weight}~$|g|$ of a shadow polyhedron~$(P,g)$ is the
  sum of the absolute values of the gleams on the regions of $P$.
\end{defi}

\begin{lemma}
\label{lem:sc-connectsum}
Shadow complexity is sub-additive under connected sum: for $M_1$,
$M_2$ two oriented 3-manifolds containing graphs $T_1$, $T_2$,
\[
\scomp(M_1\connect M_2, T_1 \cup T_2) \le \scomp(M_1,T_1) + \scomp(M_2,T_2).
\]
\end{lemma}
\begin{prf}
  Let~$P_1$ and~$P_2$ be two shadows for $(M_1,T_1)$ and $(M_2,T_2)$
  having the least number of vertices,
  and let $W_1$ and $W_2$ the corresponding 4-thickenings. To construct
  a shadow of the connect sum, let~$x_1$ and~$x_2$ by two points in
  regions of~$P_1$ and~$P_2$, respectively, and join them
  by an arc.  The polyhedron we get can be embedded as a shadow of the
  boundary connected sum of $W_1$ and $W_2$. This polyhedron is not
  simple so we modify the construction slightly: roughly speaking, we
  put our fingers at the two ends of the arc and the push $P_1$ towards
  $P_2$ along the arc until they meet in the middle along a disk.
  More precisely, identify a closed regular neighborhood of~$x_1$
  and~$x_2$, and put gleam~0 on the resulting disk region.
\end{prf}
\begin{question}
  Is shadow complexity additive under connected sum?
\end{question}
If the answer to the above question were ``yes'', a consequence
would be the following:
\begin{lemma}
  If shadow complexity is additive under connected sum, then for any
  closed 3-manifold $M$,
  \[\scomp(M)\leq 2c(M),\]
  where $c(M)$ is Matveev's complexity.
\end{lemma}

\begin{prf}
  Let $P$ be a minimal spine of $M$, i.e., a simple
  polyhedron whose 3\hyp thickening is homeomorphic to the complement~$M'$
  of a ball in~$M$ and containing the least possible number of
  vertices. Then $P$, equipped with gleam~0 on every region, is a
  shadow of $M'\times [-1,1]$, with boundary $M\connect
  \overline{M}$. Therefore $\scomp(M\connect \overline{M})\leq c(M)$ and
  the thesis follows.
\end{prf}

It is worth noting that the consequence of the above lemma is true for
all the 3-manifolds with Matveev's complexity up to $9$: we were able
to check the inequality for all of them using
Proposition~\ref{prop:arbitrarygluing} and the basic blocks exhibited
by Martelli and Petronio~\cite{MP01:3ManifoldsComplexity9}.

We next show that shadow complexity does not increase under
Dehn surgery.
\begin{lemma} \label{lem:surgery}
  Let $L$ be a framed link contained in an oriented
  3-manifold $M$ and $P$ be a shadow of $(M,L)$. A manifold $M'$
  obtained by Dehn surgery on $L$ has a shadow obtained by
  capping each component of $\partial P$ by a disk.
\end{lemma}

\begin{prf}
Let~$W$ be a 4\hyp thickening of~$P$.
Surgery of~$M$ along a component of $L$ with integer
coefficients corresponds to gluing a 2-handle to~$W$.  Gluing the core
of this 2-handle to~$P$ gives a shadow of~$W$, and the definition of
Dehn surgery on a framed link ensures that the gleam on the capped
region does not change.
\end{prf}

\begin{rem}
  Lemma~\ref{lem:surgery} together with the projection construction
  described in Section~\ref{sec:relativeshadows} give an easy
  proof that any closed 3-manifold has a shadow, since any 3-manifold can be
  presented by an integer surgery on a link in $S^3$.
\end{rem}

\begin{prop}\label{prop:arbitrarygluing}
  Let $M_1$ and $M_2$ be two oriented manifolds such that both
  $\partial M_1$ and $\partial M_2$ contain torus components~$T_1$
  and~$T_2$. Let $P_1$ and $P_2$ be shadows of $M_1$ and $M_2$, and
  let $M$ be
  any 3-manifold obtained by identifying $T_1$ and $T_2$ with an
  orientation\hyp presevering homeomorphism. Then $M$ has a shadow
  which can be obtained from $P_1$ and $P_2$ without adding any new
  vertices. In particular, any Dehn filling of a 3-manifold can be
  described without adding new vertices.
\end{prop}

\begin{prf}
  Let~$W_1$ and~$W_2$ be the 4-thickenings of~$P_1$ and~$P_2$. The
  tori~$T_1$ and~$T_2$ are equipped with the meridians~$\mu_1$
  and~$\mu_2$ of the external boundary components~$l_1$ and~$l_2$
  of~$P_1$ and~$P_2$. Also fix longitudes~$\lambda_i$ on
  them. The orientation\hyp reversing homeomorphism identifying~$T_2$
  and~$T_1$ sends~$\mu_2$ into a simple curve $a\lambda_1 +b\mu_1$
  and~$\lambda_2$ into a curve $c\lambda_1 +d\mu_1$.

  We now describe how to modify~$P_1$ and construct a shadow~$P'_1$
  of~$M_1$ embedded in a new 4-manifold~$W'_1$ such that the meridian
  induced by~$P'_1$ on~$T_1$ is the curve~$a\lambda_1+b\mu_1$. To
  construct $P'_1$ let us construct a shadow of the Dehn filling of
  $M_1$ along $T_1$ whose meridian is $a\lambda_1 +b\mu_1$.  It is a
  standard fact that any surgery on a framed knot can be translated
  into an integer surgery over a link as shown in
  Figure~\ref{fig:rationaltointeger}.

\begin{figure}
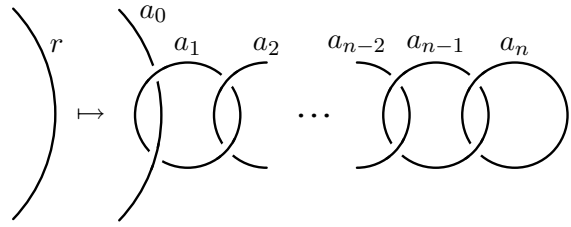


  \[
  \mfigb{draws/shadow.30} \mapsto \mfigb{draws/shadow.31}
  \]
  \caption{In this picture we show how to transform a rational surgery
    with coefficient $r$ over a knot into an integer surgery over a
    link. The coefficient $a_i$ are those of the continued fraction
    expansion of $r$, namely those of the equality:
    \[r=a_0-\frac{1}{a_1-\frac{1}{a_2\dots -\frac{1}{a_n}}}\]}
  \label{fig:rationaltointeger}
\end{figure}

\begin{figure}
  \centering
  \includegraphics{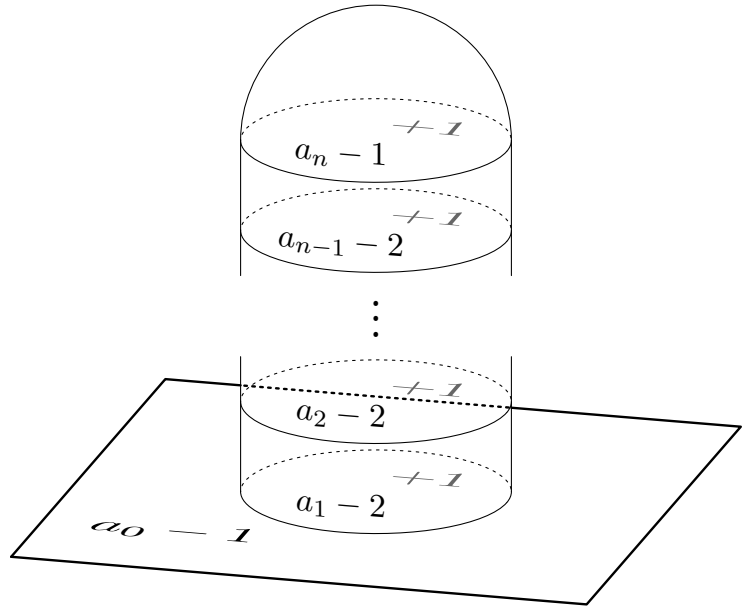}
  \caption{The shadow with no vertices corresponding to the surgeries
    on the chain in Figure~\ref{fig:rationaltointeger}.
    Intrinsically in the 4-manifold, this is equivalent to a chain of
    spheres, each intersecting the next in just one point.}
  \label{fig:rationalshadow}
\end{figure}

With the notation of the figure, we glue $n$~copies of the
polyhedron~$H$ of Example \ref{examp:hopflink} to~$P_1$ so that one component
of $\partial H_1$ is identified with $l_1$, $H_j$ is glued to
$H_{j+1}$ and $H_{j-1}$ and the free component of $\partial H_n$ is a
knot $l'_1$. On the level of the boundary of the thickening we are
gluing $n$ copies of the complement of the Hopf link in
$S^3$ to~$M$.  Since the complement of the Hopf link is $T^2\times
[0,1]$ the final 3-manifold is
unchanged. But now the polyhedron $P'_1=P_1\cup H_1\cup\ldots \cup
H_n$ can be equipped with gleams so to describe the operation of
Figure~\ref{fig:rationaltointeger}. The
meridian $\mu'_1$ of $l'_1$ is now by construction the curve which,
expressed in the initial base of $T_1$, is $a\lambda_1+b\mu_1$. Hence
we can now glue $P'_1$ and $P_2$ along $l'_1$ and $l_2$ and choose
suitably the gleam of the region of $P'_1\cup P_2$ to obtain the
desired homeomorphism.
\end{prf}

\begin{cor}
  \label{cor:sc-torus-sum}
  Shadow complexity is sub-additive under torus sums, i.e., under the
  gluing along toric boundary components through orientation reversing
  homeomorphisms.
\end{cor}

One special case of torus sums is surgery, gluing in a solid torus.
Surgery can decrease any reasonable notion of complexity, so
shadow complexity cannot be additive under torus sums in general,
hence we ask the following:
\begin{question}
  Let $M_1$ and $M_2$ be two oriented 3-manifolds with
  incompressible torus boundary components $T_1$ and $T_2$. Is it true
  that any torus sum of $M_1$ and $M_2$ along $T_1$ and $T_2$ has
  shadow complexity equal to $\scomp(M_1)+\scomp(M_2)$?
\end{question}
\subsection{Complexity zero shadows}
In this subsection we classify the
manifolds having zero shadow complexity.
\begin{defi}
  An oriented 3-manifold is said to be a {\it graph manifold} if it
  can be decomposed by cutting along tori into blocks homeomorphic to
  solid tori and $R\times S^1$, where $R$ is a pair of pants (i.e., a
  thrice\hyp punctured sphere).
\end{defi}

Graph manifolds can also be characterized as those manifolds which
have only Seifert-fibered or torus bundle pieces in their JSJ
decomposition.

\begin{prop}[Complexity zero manifolds]\label{prop:complexity-zero}
  The set of oriented 3-manifolds admitting a shadow containing no
  vertices coincides with the set of oriented graph manifolds.
\end{prop}
\begin{prf}
To see that any graph manifold has a shadow without
vertices, notice that a disk with boundary colored by~$e$
is a shadow of a solid torus and, similarly, a pair
of pants~$R$ is a shadow of $R\times
S^1$. Proposition~\ref{prop:arbitrarygluing}
shows that any gluing of these blocks can be described by a shadow
without vertices.

For the other direction, we must show that if a 3-manifold~$M$ has a
shadow~$P$
without vertices, then it is a graph manifold. The polyhedron $P$ can
be decomposed into basic blocks as follows. Since $P$ contains no
vertices, a regular neighborhood of $\Sing(P)$ in $P$ is a disjoint
union of blocks of the following three types:
\begin{enumerate}
\item the product of a $Y$-shaped graph and $S^1$;

\item the polyhedron obtained by gluing one boundary component of an
  annulus to the core of a M\"obius strip; and

\item the polyhedron obtained by considering the product of a
  $Y$-shaped graph and $[-1,1]$ and identifying the graphs $Y\times
  \{1\}$ and $Y\times \{-1\}$ by a map which rotates the legs of the
  graph of $\frac{2\pi}{3}$.
\end{enumerate}
Let $\pi:M\to P$ be the projection of $M$ on $P$.  The complement of
the above blocks in $P$ is a disjoint union of (possibly
non-orientable) compact surfaces. The preimage under $\pi$ of each of
these surfaces is a (possibly twisted) product of the surface with
$S^1$ and hence is a graph manifold. Moreover, the preimage under
$\pi$ of the above three blocks is a 3-dimensional sub-manifold of $M$
which admits a Seifert fibration (induced by the direction parallel to
$\Sing(P)$) and hence is graph manifold.
\end{prf}

\subsection{Decomposing shadows}
\label{sec:decomposing-shadows}
In Proposition~\ref{prop:complexity-zero}, we saw how to decompose a
shadow with no vertices into elementary pieces.  For more general
shadows, we will need a new type of block.  For simplicity, we will
suppose that the boundary of~$P$ is all marked ``external'' and that
the singular set $\Sing(P)$ of the shadow~$P$ is connected and
contains at least one vertex.  (This last can always be achieved by
modifying $P$ with suitable local moves.)  Let~$P$ be a shadow for a
3-manifold $M$, possibly with non empty boundary, and let $\pi:M\to P$
be the projection.  Then we have the following:
\begin{prop}\label{prop:blocks}
  The combinatorial structure of $P$ induces through $\pi^{-1}$ a
  decomposition of $M$ into blocks of the following three
  types:
\begin{enumerate}
\item products $F\times S^1$ where $F$ is an orientable surface, or
  $F\mathbin{\tilde{\times}} S^1$ with $F$ non\hyp orientable;
\item products of the form $R\times [-1,1]$, where $R$ is a pair of
  pants; and
\item genus 3 handlebodies.
\end{enumerate}
\end{prop}
\begin{prf}
Decompose the polyhedron~$P$ by taking regular neighborhoods of the
vertices and then regular neighborhoods of the edges in the
complement of the vertices.  This decomposes~$P$ into blocks of the
following three types:
\begin{enumerate}
\item surfaces (corresponding to the regions);
\item pieces homeomorphic to the product of a $Y$-shaped graph and
  $[-1,1]$; and
\item regular neighborhoods of the vertices.
\end{enumerate}

The preimage of the first of these blocks is a
block of the first type in the statement. Let us consider the preimages
of the products $Y\times [-1,1]$. The 4-dimensional
thickening of one of these blocks is the product of the 3-dimensional
thickening~$\mathbf Y$ of the $Y$-graph and $[-1,1]$, where $\mathbf Y$ is a
3-ball containing a properly embedded copy of $Y$ and collapsing on
it. The preimage in~$M$ of this block is the product of $[-1,1]$ with
$\partial {\mathbf Y} - \partial Y$, which is a pair of pants.

Let us denote by $V$ the simple polyhedron formed by a regular
neighborhood of a vertex in $P$.  We are left to show that
$\pi^{-1}(V)$ is a genus 3-handlebody.  The 4-thickening of $V$ is
${\mathbf V} \times [-1,1]$, where $\mathbf V$ is the 3-dimensional thickening
of $V$, i.e., a 3-ball into which $V$ is properly embedded. In
particular, $\partial V\subset \partial \mathbf V$ is a tetrahedral graph
and so $\partial \mathbf V$ is split into four disks by $\partial V$. One
can decompose $\partial( {\mathbf V}\times [-1,1])$ as $\partial {\mathbf V}
\times [-1,1]\cup \partial {\mathbf V}\times \{-1,1\}$. The part of this
boundary corresponding to $M$ is the complement of $\partial V$ and is
homeomorphic to ${\mathbf V}\times \{-1\}\cup (\partial {\mathbf V}-
\partial V)\times [-1,1] \cup {\mathbf V}\times \{1\}$; this is composed
of two 3-balls connected through 4 handles of index one (each of which
corresponds to one of the four disks into which $\partial V$ splits
$\partial \mathbf V$).
\end{prf}

Note that boundary the blocks of the second two types in
Proposition~\ref{prop:blocks} are themselves naturally decomposed into
annuli and pairs of pants.

\subsection{A family of universal links}
\label{sec:family-univ-links}
Now suppose further that~$P$ (and therefore~$M$) has no boundary,
and consider the union of the blocks of the second two types in
Proposition~\ref{prop:blocks}.  These two types of blocks meet in
pairs of pants, and the remaining boundary is obtained from the
annuli; therefore, we are left with a manifold~$S_P$ with
boundary a union of tori, which depends only on the polyhedron~$P$ and
not on the gleams.  (The original manifold~$M$ can be obtained by
surgery on~$S_P$.)

In this subsection we show that~$S_P$ is a hyperbolic cusped 3-manifold
whose geometrical structure
can be easily deduced from the combinatorics of $P$. We furthermore
show how to present $S_P$ as the complement of a link in a connected
sum of copies of $S^2\times S^1$.
 
As before, let $P$ be a simple polyhedron (now with no boundary) such
that $\Sing(P)$ is connected and contains at least one vertex; let
$c(P)$ be the number of vertices. Let
$S(P)$ be the regular neighborhood
of $\Sing(P)$ in~$P$, which we think of as a simple polyhedron with
boundary colored ``internal''.  Let $l_1,\ldots,l_k$ be the components of
$\partial S(P)$ in $P$. To each $l_i$ we assign a positive integer
number $c_i$ called its {\it valence} by counting the number of
vertices touched by the region $R_i$ of $S(P)$ containing $l_i$ and an
element of $\mz_2$ given by the $\mz_2$-gleam $g_i$ of the region of
$S(P)$ containing $l_i$.

Let $X_P$ be the 4-thickening of $S(P)$ provided by Turaev's
Reconstruction Theorem; $X_P$ collapses onto a graph with Euler
characteristic $\chi (S(P))=-c(P)$ and so $\partial X_P$ is a
connected sum of $c(P)+1$ copies of $S^2\times S^1$. Moreover,
$\partial S(P)$ is a link $L_P$ in $\partial X_P$.  The manifold~$S_P$
introduced earlier is the complement of~$L_P$ in~$\partial X_P$.

$S_P$ has a natural hyperbolic structure which we can understand in
detail, as we will now see.

\begin{prop}
  \label{prop:hyperbolicpieces}
  For any standard shadow surface~$P$,
  $S_P$ can be equipped with a complete, hyperbolic metric with
  volume equal to $2v_{\oct} c$.
\end{prop}

\begin{prf}
The main point of the proof is to construct an hyperbolic
structure on a block corresponding to a vertex in $S(P)$ and then to
show that these blocks can be glued by isometries along the edges
of~$S(P)$.

Let us realize a block of type 3 as follows. In $S^3$, pick two
disjoint 3-balls $B_0$ and $B_\infty$ forming neighborhoods
respectively of 0 and $\infty$. Connect them using four 1-handles
$L_i$, $i=1,\ldots,4$, positioned symmetrically, as shown in
Figure~\ref{fig:hyperbolicstructure}.

In the boundary of the so obtained genus 3-handlebody consider the 4
thrice\hyp punctured spheres formed by regular neighborhoods of the
theta-curves connecting $B_0$ and $B_\infty$ each of which is formed
by 3-segments parallel to the cores of three of the 1-handles. These
four pants are the surfaces onto which the blocks of type 2 in
Proposition~\ref{prop:blocks} are to be glued. Indeed, these blocks
are of the form
$R\times [-1,1]$ where $R$ is a thrice punctured sphere, and they are
glued to the blocks of type~3 along $R\times \{-1,1 \}$. We will now
exhibit an hyperbolic structure on this block so that these 4
thrice\hyp punctured spheres become totally geodesic and their complement
is formed by $6$ annuli which are cusps of the structure.

Consider a regular tetrahedron in $B_0$ whose barycenter is the center
of $B_0$ and whose vertices are directed in the four directions of the
1-handles $L_i$. Truncate this tetrahedron at its midpoints as shown
in Figure~\ref{fig:hyperbolicstructure}. The
result is a regular octahedron~$O_0$ contained in~$B_0$, with 4~faces
(called ``internal'') corresponding to the vertices of
the initial tetrahedron and 4~faces (called ``external'') corresponding
to the faces of the initial tetrahedron. Do the same construction
around $\infty$ and call the result $O_\infty$. The handlebody
$B_0\cup B_\infty \cup L_i$, $i=1,\ldots,4$ can be obtained by gluing
the internal faces of $O_0$ to the corresponding internal faces of
$O_\infty$.  The remaining parts of the boundaries of the two
octahedra are four spheres each with three ideal points and
triangulated by two triangles. If we put the hyperbolic
structure of the regular ideal octahedron on both~$O_0$ and~$O_\infty$,
then, after truncating with horospheres near the vertices, we get the
hyperbolic structure we were searching for: the geodesic thrice\hyp
punctured spheres come from the boundary spheres without their cone
points and the annuli are the cusps of the structure. Each (annular)
cusp has an aspect ratio of $\frac{1}{2}$ since it is the union of two
squares, the sections of the cusps of an ideal octahedron near a
vertex.
\begin{figure}
\psfrag{B_0}{$B_0$}
\psfrag{B_infty}{$B_\infty$}
\psfrag{l_1}{$L_1$}
\psfrag{l_2}{$L_2$}
\psfrag{l_3}{$L_3$}
\psfrag{l_4}{$L_4$}
  \centerline{\includegraphics[width=7.4cm]{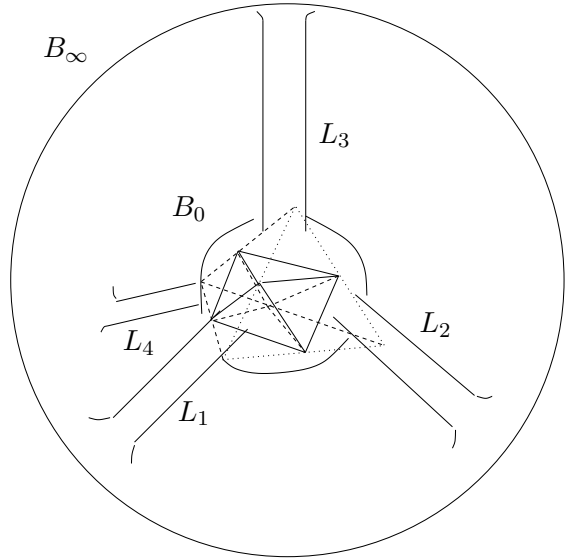}}
  \caption{In this picture we show how to connect the two balls $B_0$
    and $B_\infty$ in $S^3$ using the four legs $L_i$, $i=1,\ldots,4$.
    In the center of $B_0$, we visualize how the regular octahedron
    $O_0$ is embedded. In the figure, the four internal faces are
    directed towards the four legs of the handlebody since they are
    identified with the four internal faces of the octahedron
    $O_\infty$.}
  \label{fig:hyperbolicstructure}
\end{figure}
To show that these blocks can be glued and form a hyperbolic
manifold $S_P$ it suffices to notice that the thrice\hyp punctured
spheres in a block of type 3 are all isometric.
\end{prf}

\begin{prop}
  The Euclidean structure on the cusp corresponding to a boundary
  component~$l_i$ of~$S(P)$ is
  the quotient of
  $\mr^2$ under the two transformations $(x,y)\cong (x+2,y)$ and,
  $(x,y)\cong (x+g_i,y+c_i)$.
\end{prop}

\begin{prf}
The cusp corresponding to $l_i$ is obtained by gluing some
of the annular cusps in the blocks of the vertices: each time $l_i$
passes near a vertex $v$ of $S(P)$, we glue the annular cusp
corresponding to $l_i$ in the block of $v$ (note indeed that in this
block there are exactly $6$ cusps, one for each of the six regions
passing near the vertex). Since each annular cusp has a section which
is an annulus whose core has length 2 and height is 1, following $l_i$
and gluing the cusps corresponding to the vertices we meet, we
construct an enlarging annular cusp; when we conclude a loop around
$l_i$, we glued $c_i$ cusps and we got an annulus whose core has
length 2 and whose height is $c_i$. Then, we are left to glue the two
boundary components of this annulus to each other, and the
combinatorics of $S(P)$ forces us to do that by applying $g_i$ half
twists to one of the two components.
\end{prf}

Finally, we give a more explicit description of the link~$L_P$ in
terms of surgery on $S^3$.

\begin{figure}[ht!]
\psfrag{0}{$0$}
  \centerline{\includegraphics[width=10.4cm]{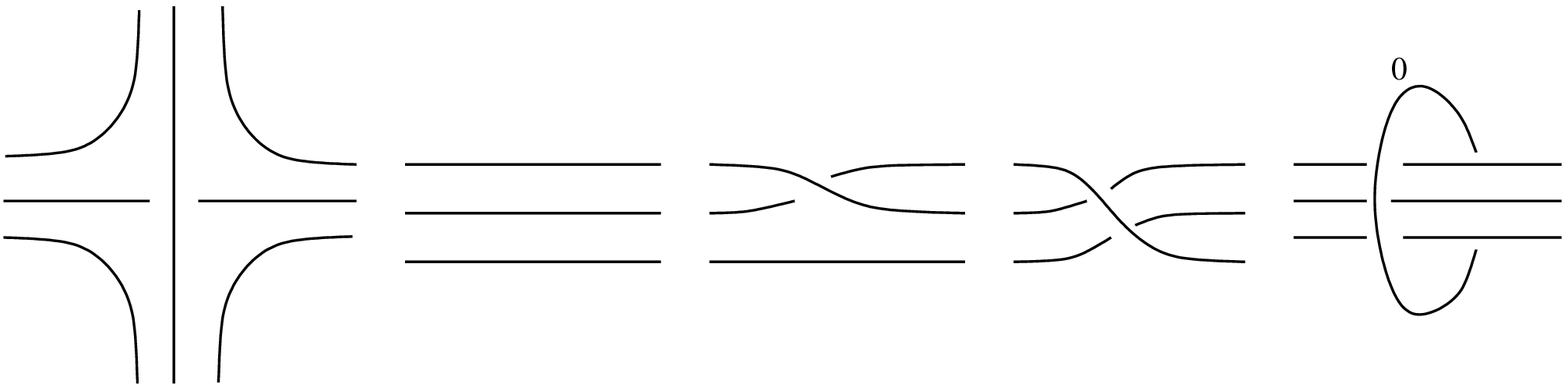}}
  \centerline{\includegraphics[width=8.4cm]{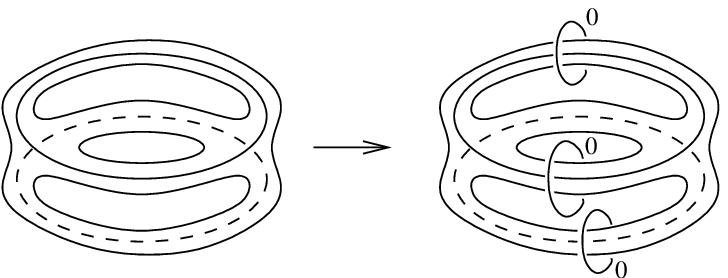}}
  \caption{In the upper part of the picture we draw the basic blocks
    composing the links $L_P$. In the lower part we work out an
    example: at the left part we show an example of $S(P)$ (the curves
    represent $\partial S(P)$). In the right part we encircle the
    complement of a maximal tree of $\Sing(P)$ with 0-framed
    meridians. }
  \label{fig:blocklink}
\end{figure}

\begin{prop}
  \label{prop:universallink}
  $S_P$ can be presented as the complement of a link~$L_P$ in the
  manifold obtained by surgering $S^3$ over a set of $c(P)+1$
  unknotted 0-framed meridians (where $c(P)$ is the number of vertices
  of $P$). Moreover this link can be decomposed into blocks like those
  shown in Figure~\ref{fig:blocklink}.
\end{prop}

\begin{prf}
Let $T$ be a maximal tree in $S(P)$ and consider its
regular neighborhood, a contractible
sub-polyhedron $P'$ of $S(P)$; $S(P)$ can be
recovered from $P'$ by gluing to $P'$ the blocks corresponding to the
edges of $S(P)\setminus T$.  Let $B$ be the 3\hyp dimensional
thickening of~$P'$, and let $X'$ be $B\times[-1,1]$, the 4\hyp
dimensional thickening of~$P'$.  The trivalent
graph $\partial P'$ is contained in $\partial B\times \{0\}\subset
\partial X'=S^3$.  Moreover we can push $P'$ into $\partial X'$ by an
isotopy keeping its boundary fixed. Then $\partial P'$ is the
boundary of a contractible polyhedron in $B^3$ and
hence is composed by joining some copies of the blocks shown in
the upper-left part of Figure \ref{fig:blocklink} by means of
triples of parallel strands. Each time we glue back to
$P'\subset X'$ a block corresponding to an edge of $S(P)\setminus T$, we
are gluing to $X'$ a 1-handle connecting neighborhoods of two
vertices, say $v_1$ and $v_2$, of $\partial P'$. The boundary of
the polyhedron we get that way is obtained from $\partial P'$
by connecting the strands around $v_1$ and those around $v_2$
according to the combinatorics of $S(P)$ and letting it pass over
the 1-handle once: this can be represented by a passage
through a 0-framed meridian. Performing this
construction on all the edges of $S(P)-T$ one gets the
link $L_P$ of the form described in the statement.
\end{prf}

As already noticed, any closed oriented 3-manifold has a shadow;
moreover, up to applying some basic transformations to such a shadow,
we can always suppose it to be standard. This has the following
consequence:
\begin{prop}
  The family of links $L_P$, with $P$ ranging over all standard
  polyhedra, is ``universal'': any closed orientable 3-manifold can be
  obtained by a suitable integral surgery over an element of this
  family.
\end{prop}
Since the number of standard polyhedra with at most $c$ vertices is
finite, the family of universal links $L_P$ has a natural finite
stratification given by the complexity of the polyhedron from which
each element of the family is constructed.  Using Jeff Weeks' program
SnapPea,
we were able to check that all but 4 manifolds of the cusped census
can be obtained by surgering over links corresponding to polyhedra
with at most 2 vertices.

The following inequality from the introduction is a corollary of
Gromov's results~\cite{Gromov82:Volume}:

\begin{teo}\label{teo:gromov-norm-lower}
  A 3-manifold~$M$, with boundary empty or a union of tori, has shadow
  complexity of at least $(v_{\tet}/2v_{\oct})\|M\|$.
\end{teo}

\begin{prf}
In any shadow for~$M$ with $n$ vertices, the preimage of a neighborhood
of the singular set is the disjoint union of pieces which either have the
hyperbolic structure described above (if there is at least one vertex
in the connected component) or are graph manifolds (as in
Proposition~\ref{prop:complexity-zero}).  The total Gromov
norm of these pieces is therefore $(2v_{\oct}/v_{\tet})n$.  $M$ can be
obtained from these pieces by gluing some additional pieces from the
regions: each region contributes a surface cross $S^1$.
Since the Gromov norm is non-increasing under gluing along torus
boundaries~\cite{Gromov82:Volume}, $n$ must be at least
$v_{\tet}/2v_{\oct}\cdot\|M\|$.
\end{prf}

\begin{teo}\label{teo:singularities-lower}
  A 3-manifold~$M$ with Gromov norm~$G$ has at least $G/10$ crossing
  singularities in any smooth, stable map $\pi:M\to\RR^2$.
\end{teo}
\begin{prf}
Applying the construction underlying the proof of Theorem
\ref{teo:main}, one can construct $M$ as a Dehn filling of a link
$L_P$ for a suitable simple polyhedron $P$; moreover, each singularity
of the second type (as in Figure~\ref{fig:case2}) produces a pattern
which can be
triangulated with $10$ regular ideal hyperbolic tetrahedra, and each
singularity of the first type (as in Figure~\ref{fig:case1}) can be
obtained as the
union of two regular ideal octahedra.  Hence $M$ is the Dehn filling
of an hyperbolic cusped 3-manifold whose volume is no more than
$10s{v_{\tet}}$, where $s$ is the number of crossing singularities
of~$\pi$ and ${v_{\tet}}$ is the volume of the regular ideal
tetrahedron.  The statement follows.
\end{prf}

\section{Shadows from triangulations}
\label{sec:shadows-triangulations}

In this section, we exhibit a construction which, given a
3-manifold~$M$ triangulated with $t$~tetrahedra (possibly with some
ideal vertices), produces a shadow
of the manifold containing a number of vertices bounded from
above by $kt^2$ where $k$ is a constant which does not depend on $M$.
This produces a 4-manifold
whose shadow complexity (the least number of vertices of a
shadow of the manifold) can be bounded by $kt^2$ and whose
boundary is the given 3-manifold.  Furthermore, we bound the gleam
weight and the number of 4-simplices needed to construct the
4-manifold.

Because this is the central point of the paper, we go into some
details and give an explicit estimate for $k$ and the bounds on the
gleam weights.

From now on, by a \emph{triangulation} we mean a
$\Delta$-triangulation, an assembly of simplices glued along their
faces, possibly with self-gluings.

\begin{defi}
Let $M$ be an oriented 3-manifold whose boundary does not contain
spherical components. A \emph{partially ideal}
triangulation of~$M$ is a triangulation of the singular manifold
$M/\partial M$ (obtained by identifying each boundary component to
a point) whose vertices contain the singular points corresponding
to the boundary components of~$M$.  An \emph{edge-distinct}
triangulation is a triangulation where the two vertices of each edge
(simplex of dimension~1) are different.
\end{defi}
This section is devoted to proving the following theorem which is the
main tool in proving the results announced in the introduction:
\begin{teo}\label{teo:main}
  Let $M$ is an oriented 3-manifold, possibly with boundary,
  and let $T$ be a
  partially ideal, edge-distinct triangulation of $M$ containing $t$
  tetrahedra. There exists a shadow $P$ of $M$ 
  which is a boundary\hyp decorated standard polyhedron, contains at most $18t^2$
  vertices, and has gleam weight at most $108t^2$.
\end{teo}

\begin{cor}
  With the assumptions an in Theorem~\ref{teo:main}, except that $T$
  is not necessarily edge-distinct, then $M$ has a shadow which is a
  standard polyhedron, contains at most $24^2\cdot18t^2$ vertices, and
  has gleam weight at most $24^2\cdot 108t^2$.
\end{cor}
\begin{prf}
  Apply Theorem~\ref{teo:main} to the barycentric subdivision of~$T$,
  which has $24t$~tetrahedra and is edge-distinct.
\end{prf}

\begin{prf}[Proof of Theorem~\ref{teo:main}]
The main idea of the proof is to pick a map from
$M/\partial M$ to $\mr^2$, stabilize its singularities and
associate to this map its Stein factorization, which turns out to be a
decorated shadow of $M$. We split the proof of the theorem into
6~main steps.
\begin{enumerate}
\item Define an initial projection map.  We map all the vertices to
  the boundary of the unit disk, so that they don't interfere with the
  bulk of the construction.

\item Modify the projection map to get a map which is stable in the
  smooth sense.  This involves modifying the projection in a
  neighborhood of the edges, in a uniform way along the edge.

\item Construct the shadow surface over the complement of a 
  neighborhood of the codimension~2 singularities of the map.  Here
  the shadow surface is just the Stein factorization as described in
  the introduction.

\item Extend the construction to the neighborhoods of the
  codimension 2 singularities.  This involves analyzing the two
  interesting types of singularities.  For one of the singularities we
  modify the Stein factorization slightly to get a shadow surface.

\item Estimate the complexity of the resulting shadow.  Essentially,
  the vertices may come from interactions between a pair of edges, and
  there are quadratically many such interactions.

\item Estimate the gleams on the regions of the shadow.

\end{enumerate}

\subsection{The initial projection}
\label{sec:project}

Pick a generic map~$\pi$ from the vertices
$v_1,\ldots,v_n$ of $T$ to the unit circle in $\mr^2$ and call
$p_1,\ldots,p_n$ their images. Extend~$\pi$ to all of~$T$
in a piecewise-linear fashion to a map from $M/\partial M$ to the unit
disk (see Figure~\ref{fig:projection1}). Pick a small disk around
each~$p_i$, and let $M'$ be the complement in
$M/\partial M$ of the preimage of these disks; $M'$ is homeomorphic to $M$
minus a ball around each non-ideal vertex.

Let~$G$ be the image (via~$\pi$) in $\mr^2$ of the union of the edges
of~$T$. If $p$ is any point in $\mr^2\setminus G$ then $\pi^{-1}(p)$ is a
set (possibly empty) of circles in~$M'$ since it is a union of
segments properly embedded in the tetrahedra of $T$ never meeting the
edges of $T$.  We may think that the set of ``critical values''
of~$\pi$ is contained in~$G$.

\begin{figure}
\psfrag{v_1}{$v_1$}
\psfrag{v_2}{$v_2$}
\psfrag{v_3}{$v_3$}
\psfrag{v_4}{$v_4$}
\psfrag{p_1}{$p_1$}
\psfrag{p_2}{$p_2$}
\psfrag{p_3}{$p_3$}
\psfrag{p_4}{$p_4$}

  \centerline{\includegraphics[width=6.4cm]{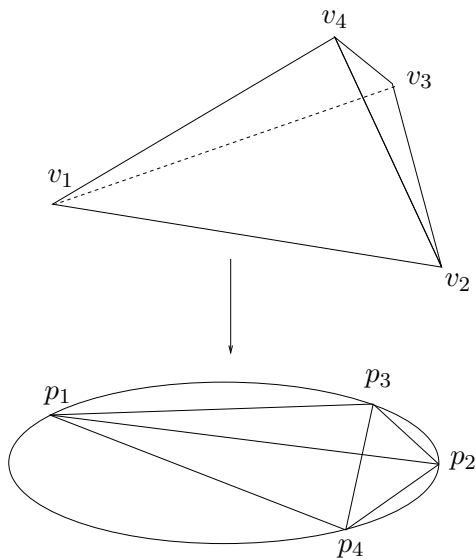}}
  \caption{In this figure we summarize the notation we fixed in the
    first step. The tetrahedron spanned from the vertices $v_i$ is a
    tetrahedron of the triangulation $T$ of $M/\partial M$. }
  \label{fig:projection1}
\end{figure}

\subsection{A stable projection}
\label{sec:stable}

The boundary of $M'$ in $M/\partial M$ is a union of ``vertical''
surfaces, surfaces which project to segments of small circles around the
points~$p_i$. In the next sections
we will restrict ourselves to $M'$ and construct a Stein
factorization of $\pi:M'\to \mathbb{R}^2$, which turns out to be a
shadow of $M'$; we will then modify it to get a shadow for~$M$.

The image of the edges~$e_i$, $i=1,\ldots, r$ of~$T$ form a set of
segments~$f_i$ in the unit circle. Since two edges~$e_i$ and~$e_j$
could have the same endpoints in $T$, some $f_i$ could coincide. To
avoid this, we modify~$\pi$ slightly around small regular
neighborhoods of the edges in $M'$ so that the projections of
different edges with the same endpoints in~$T$ are distinct segments in
the unit circle running parallel to each other. This can be done by
operating in disjoint small cylindrical neighborhoods of the edges,
since no vertices of~$T$ are contained in~$M'$.
The resulting map is no longer PL.

Let us keep calling~$G$ the graph which is the (modified) image of
the~$f_i$.
The edges of~$G$ are now straight segments away from neighborhoods of
the~$p_i$, with bends near the~$p_i$.  We now study the behavior of
the projection map on a cylindrical regular
neighborhood~$C_i$ of each edge~$e_i$ of $T$.  Transverse to~$e_i$
in~$C_i$ is a triangulated disk with one interior vertex from $e_i$
and triangles coming from the tetrahedra of $T$ incident to $e_i$;
the projection of this disk in $\mathbb{R}^2$ is a segment transverse
to $f_i$.  Let~$\pi_t$ be the map from this transverse disk to the
transverse segment.  The map from $C_i$ to the neighborhood of $f_i$
is the product of~$\pi_t$ with an interval, hence it suffices to
study~$\pi_t$.

For instance, consider the following possibility for~$\pi_t$: let $Q$
be a square
triangulated into four triangles by coning from the center, let
$A$, $B$, $C$, $D$ be its vertices in cyclic order, and consider the
PL map from~$Q$ to $[-1,1]$
sending~$A$ and~$C$ to~$1$, $B$ and~$D$ to~$-1$ and the center to~$0$. The
preimage of a point near~$1$ (resp.~$-1$) is a pair of segments near~$A$
and~$C$ (resp.~$B$ and~$D$).  The preimage of~$0$ is the cone from the
center of~$Q$ to the midpoints of its edges. This map is the typical
example of a saddle on the base of $C_i$.

If we repeat the above
construction with an hexagon, sending the vertices alternately to~1
and~$-1$, we obtain instead a ``monkey saddle'', which is not stable
(from the smooth point of view). The inverse image of 0 is a cone over
the midpoints of the edges from the center: a six-valent
star. As shown in Figure~\ref{fig:stabilization}, in this
case $\pi_t$ can be perturbed to a map having two stable critical
points as shown in the figure.
In general, if the inverse image of
a critical value is a star with~$2k$ legs, then~$\pi_t$ can be
perturbed to a map containing $k-1$ stable saddle points all having
distinct images in the segment.

There is one case left: when the
whole disk is projected on one side of 0 in $[-1,1]$.  In this case the
singular point in the center of the disk is an extremum and we keep
the map unchanged.

\begin{figure}
\psfrag{Preimage of zero}{Preimage of zero}
\psfrag{Preimage of -1/2}{Preimage of $-1/2$}
\psfrag{Preimage of 1/2}{Preimage of $1/2$}
\psfrag{1}{$1$}
\psfrag{-1}{$-1$}
  \centerline{\includegraphics[width=7.4cm]{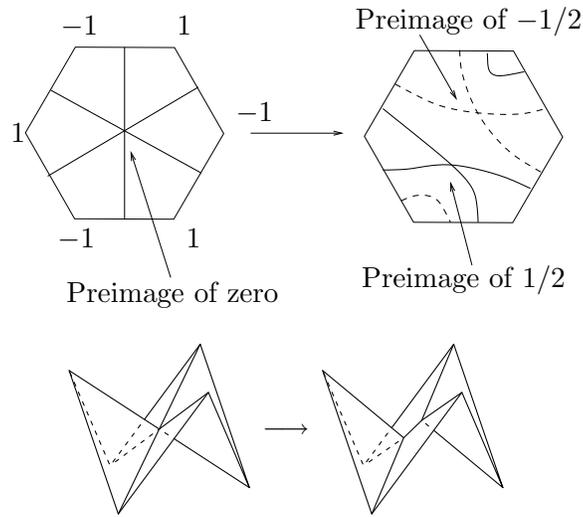}}
\vspace{10pt plus 3pt minus 3pt}
  \centerline{$\mfigb{draws/surface.80}\longrightarrow\mfigb{draws/surface.81}$}
  \caption{Stabilizing the map sending the vertices of an
  hexagon over $[-1,1]$ alternatively in $\pm 1$. In the right part
  we show the preimages of the two singular values ($\pm\frac{1}{2}$) of
  the stabilized map.  The lower diagram is another representation,
  where the projection to~$\mr$ is the
  projection onto the vertical axis.}
  \label{fig:stabilization}
\end{figure}

We now modify $\pi:C_i\to \mathbb{R}^2$ as above around each edge
$e_i$ to get the cylinder of a stable map from a disk to a
segment. This increases the set of critical values of $\pi$ near $f_i$
so that it no longer coincides with $f_i$, but is formed by a set of
strands running parallel to it, all corresponding to stable
singularities of the map. Let us keep calling $G$ the graph in $\mr^2$
made of these critical values; as before, it consists of straight
segments away from the vertices~$p_i$, with bent segments near the~$p_i$.
These straight segments are cut by their intersections $x_k,\ k=1,\ldots, l$ 
(which, together with the points~$p_i$, form the vertices of $G$) into
sub-segments $h_j,j=1,\ldots, m$ (which form the edges of $G$).

\subsection{The Stein factorization away from codimension~2 singularities}
\label{sec:construct-shadow-1}
Pick a regular
neighborhood of each vertex of~$G$, and let $M''$ be the 
preimage through $\pi$ of the complement of these neighborhoods.
We will now construct a shadow~$P''$ for~$M''$ from the Stein
factorization for the map~$\pi$, as shown in
Figure~\ref{fig:projection2}.

Let $R_1$, \dots, $R_m$, $R_{\infty}$ be the connected components of
$\mr^2\setminus G$, where $R_{\infty}$ is the unbounded region and $R_i$,
$i\neq \infty$, are disks.
By construction, $\pi^{-1}(R_{\infty})$ is
empty and $\pi^{-1}(R_i)$, $i\neq \infty$ is a
disjoint union of $n_i$ open solid tori in $M'$.  
For an edge $h_j$ of $G$,
let $\alpha_j$ be a small arc intersecting it
transversally and connecting two regions,
say~$R_0$ and~$R_1$. Let $q_0$ and
$q_1$ be the endpoints of $\alpha_j$; the (possibly disconnected)
surface $S_j=\pi^{-1}(\alpha_j)\subset M''$ is the cobordism between
$\pi^{-1}(q_0)$ and $\pi^{-1}(q_1)$ whose possible shapes are depicted
in Figure~\ref{fig:projection2}.

To construct the Stein factorization of $\pi$, for each $R_i$, take
$n_i$~copies of~$R_i$.
We need to connect these regions to each other near the
centers of the segments~$h_j$. To do this, we apply the procedure of
Figure~\ref{fig:projection2}, where all the possible behaviors of $S_j$
are examined. Saddle singularities produce a
singular set in the polyhedra used to connect the regions.
Shrinking singularities (when the transverse map $\pi_t$ in the
previous step maps entirely on one side of the singularity) produce a
boundary segment of~$P''$, which we mark as ``false''; temporarily mark
the rest of the boundary of~$P''$ as ``internal''.
Call the regions which are involved in the singularity or boundary
over~$h_j$ the \emph{interacting regions}.
\begin{figure}
\psfrag{R_1}{$R_1$}
\psfrag{R_2}{$R_2$}
\psfrag{f_1}{$f_1$}
\psfrag{f_2}{$f_2$}
\psfrag{Singular values}{Singular values}
\psfrag{P}{$P$}
\psfrag{Sing(P)}{$\Sing(P)$}
\psfrag{pigreco'}{$\pi_1$}
\psfrag{pigrecop}{$\pi_2$}
\psfrag{R_1^i}{$R_1^i$}
\psfrag{R_1^j}{$R_1^j$}
\psfrag{R_1^k}{$R_1^k$}
\psfrag{R_1^l}{$R_1^l$}
\psfrag{R_2^i}{$R_2^i$}
\psfrag{R_2^j}{$R_2^j$}
  \centerline{\includegraphics[width=12.4cm]{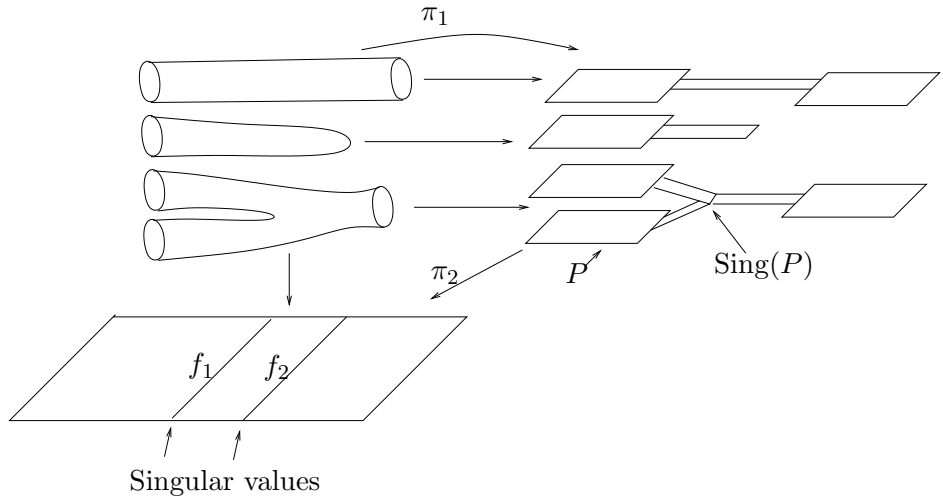}}
  \caption{The construction we perform in
    Step~3. To each fiber over a region of $\mr^2$ we associate a
    region of~$P''$. We glue these regions along their boundaries
    to the type of singularity which intervenes. In the highest part
    of the picture we see the
    simplest case: when the singularity does not affect the fibers
    corresponding to two regions of~$P''$. The second case is the one of
    a shrinking singularity, which creates a false boundary component
    in~$P''$.  The last case is the case of a simple saddle singularity,
    which creates an arc in $\Sing(P'')$. }
  \label{fig:projection2}
\end{figure}

If we repeat the above construction for all pairs of regions in
contact through a segment of the family $h_j$, we get a decorated simple
polyhedron which represents the Stein factorization of
$\pi:M''\to \mathbb{R}^2$.  We can naturally find
maps $\pi_1:M''\to P''$ and $\pi_2:P''\to \mr^2$ so that $\pi
=\pi_2 \circ \pi_1$.

Let us analyze $\partial_e P''$.  Currently, $\pi_2(P'')$
covers the complement in $\mr^2$ of small circular neighborhoods
of the vertices of~$G$, which are either the points~$p_i$ (the images
of the vertices of~$T$) or
intersections~$x_k$ of the edges~$f_i$.
The inverse image in~$P''$ of the boundaries of these circular
neighborhoods is~$\partial_e P''$, which is a trivalent graph possibly
with some free ends.

\subsection{Codimension-2 singularities}
\label{sec:extend-shadow}

We now describe how to extend~$P''$ to get a shadow~$P'$ of~$M'$.
We fill in the gaps of the polyhedron near the
intersections~$x_k$ of critical values by using a simple polyhedron with at
most 2 vertices per intersection.
 
Near each~$x_k$ two segments of critical values
intersect, say~$s_1$ and~$s_2$.
For every region which is not interacting over either~$s_1$ or~$s_2$,
we fill in the hole over $x_k$ with a disk.  Also, if the region(s)
interacting over~$s_1$ do not meet those interacting over~$s_2$, we
can fill in the hole with the same simple blocks as in
Figure~\ref{fig:projection2} without adding any new vertices.  In
particular this always occurs if the singularity over~$s_1$ or~$s_2$
is a shrinking singularity.

We are left with the case when both~$s_1$ and~$s_2$ correspond to
saddle singularities.  In this case the component of the
preimage~$\pi^{-1}(x_k)$ containing the singularities is a
connected 4-valent graph~$F$ in $M'$ with two
vertices, one from each singular segment.  The edges of~$F$ are
oriented: At a generic point on
an edge of~$F$, a small disk in~$M'$ transverse
to the edge maps homeomorphically to its image
in~$\mr^2$ and so we can pull back the orientation of~$\mr^2$
to it. Then, since $M'$ is oriented,
we can orient the edges of~$F$. Each vertex of $F$ corresponds to a
codimension 1 singularity whose singular values are contained either
in $s_1$ or in $s_2$. Moreover, near each vertex of $F$, two edges are
incoming and the other two outgoing. Therefore the only
possibilities for~$F$ are these graphs:
\[
\mfigb{draws/singularity.5}\textrm{\quad or\quad}\mfigb{draws/singularity.51}.
\]
While passing through the codimension 1 singularity corresponding to a
vertex, the edges of~$F$ recouple so that incoming edges are glued
to outgoing edges.

We now analyze these two cases and show that if $F=\pi^{-1}(x_k)$ has the
first shape, then its neighborhood in $M'$ can be reconstructed by
using a shadow polyhedron with one vertex, while in the second
case, two vertices are sufficient.
In both cases, the regular neighborhood of $F$ in $M'$ is a
3-handlebody~$H(F)$. By construction, the
boundary~$\Sigma_3$ of this handlebody projects,
through~$\pi_1$, to a component~$G$ of~$\partial_i P''$ and,
through~$\pi$, to a
circle in $\mr^2$ which is the boundary of a small regular
neighborhood of~$x_k$.
The thickening of~$G$, contained in the thickening of~$P''$
constructed so far, is
another 3-handlebody $H(G)$ lying vertically (through
$\pi_2$) over this circle, and whose boundary is identified in $M'$
with $\Sigma_3$.
We will show that in both cases the union of these two handlebodies
is~$S^3$ and then construct a shadow of~$(S^3,G)$, where $G$~is
considered as a subset of $H(G) \subset S^3$.

\case{1}{$F$ is~$\mfigb{draws/singularity.9}$.}
\newsavebox{\mygraphics}
\savebox{\mygraphics}{$\mfigb{draws/singularity.9}$}
\begin{figure}
  \centerline{\includegraphics{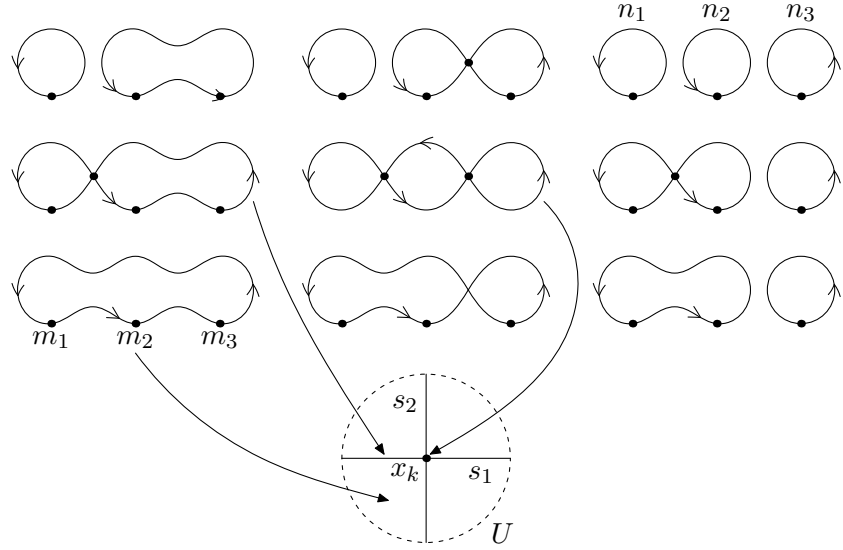}}
  \caption{\label{fig:case1}The behavior over a neighborhood~$U$ of
    the intersection of
    two segments ~$s_1$, $s_2$ of singular values of the projection
    map when the singular graph~$F=\pi^{-1}(x_k)$
    is~\usebox{\mygraphics}.  We see the different topological types
    of preimages of points, with the induced orientation.  The
    different types of points are the four regular areas (the
    components of $U\setminus(s_1\cup s_2)$), the segments $s_1$ and
    $s_2$ themselves, and the intersection $s_1\cap s_2$.  We have
    marked two different systems of curves in the fibers.  One (the
    $n_i$) are copies of the fibers over the upper right region; these
    are meridians for~$H(G)$.  The other (the $m_i$) map
    surjectively to~$\partial U$, and are meridians for~$H(F)$.  In
    the picture they appear as a choice of one point in
    each fiber.}
  \label{fig:firstsingularity}
\end{figure}
In Figure \ref{fig:firstsingularity}, we show the preimages in~$M'$ of
various points in a small circular neighborhood~$U$ of
$x_k$ in $\mr^2$. In each component of $U\setminus
(s_1\cup s_2)$ the
preimage of a point is a union of circles in $M'$. While crossing
$s_1$ or $s_2$ two arcs of these circles approach each other and after
passing through a singular position, they recouple. In the figure, we
show the preimage of a point in each of the regular areas and on each
of the singularities.

In the figure we have picked three of the edges of~$F$.  Transverse to
these edges are three meridian disks of~$H(F)$
bounded by curves $m_1$, $m_2$ and $m_3$.  Each disk can be chosen to be
a section of $\pi$ as shown by the dots in
Figure~\ref{fig:firstsingularity}, so
that each of the~$m_i$ projects homeomorphically to $\partial U$.

We now identify meridians of $H(G)$. By construction,
each circle in $M'$ which is in the preimage of a regular point in $\mr^2$
bounds a disk in~$H(G)$: the preimage of the corresponding point in the
thickening of~$P''$. Hence the meridian curves of $H(G)$ include the
circles drawn in Figure~\ref{fig:firstsingularity} over the 4 areas
near $x_k$.  In the upper-right region the preimage of a point
is composed of three circles $n_1,n_2$ and $n_3$, which we can choose
as our Heegaard system.

The handlebodies $H(F)$ and $H(G)$ are glued along their
boundaries. Since each meridian $n_i$ of $H(G)$ intersects
exactly one time one of the meridians $m_j$ of $H(F)$, we have
$H(G)\cup H(F)=S^3$.  Inside this 3-sphere, $G$ is a
3-valent graph.
The vertex of Figure~\ref{fig:singularityinspine}, with boundary
colored~$i$, forms a shadow of the pair $(S^3, G)$: Its
thickening is $B^4$ and its boundary sits in $S^3=\partial B^4$. To
see that this is a correct shadow, in
Figure~\ref{fig:firstsingembedded}, we exhibit two oriented graphs
embedded in $S^3$ which are homeomorphic respectively to $G$ and
$F$. Moreover, the graphs are equipped with systems of
meridians~$n_i$ and~$m_i$, respectively; when the meridians are
homotoped into a common
surface, the intersections $m_i\cap n_j$ coincide with the
corresponding intersections in $M'$ between the meridians of $H(F)$
and $H(G)$.

This shows that, when the codimension two singularity is
$\mfigb{draws/singularity.9}$, it is sufficient to form~$P'$ by gluing
one vertex to~$P''$ in order to extend the description of $M'$
over~$x_k$.
\begin{figure}
  \centerline{\includegraphics{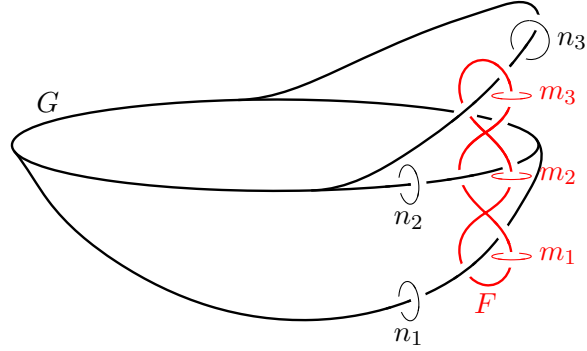}}
  \caption{The two graphs~$G$ (in black) and~$F$ (in red) embedded in
    $S^3$.  Their thickenings form two genus 3 handlebodies in $S^3$
    equipped with meridian curves whose intersections agree with
    the intersections of the curves $m_i$ and $n_i$ described in the
    text.}
  \label{fig:firstsingembedded}
\end{figure}

\case{2}{$F$ is~$\mfigb{draws/singularity.59}$.}
For this case,
Figure~\ref{fig:secondsingularity} shows the preimages in~$M'$ of
various points in a small neighborhood of~$x_k$. As shown in the
figure, two opposite areas are covered by 2~circles and the other two
by 1~circle.

\savebox{\mygraphics}{$\mfigb{draws/singularity.59}$}
\begin{figure}
  \centerline{\includegraphics{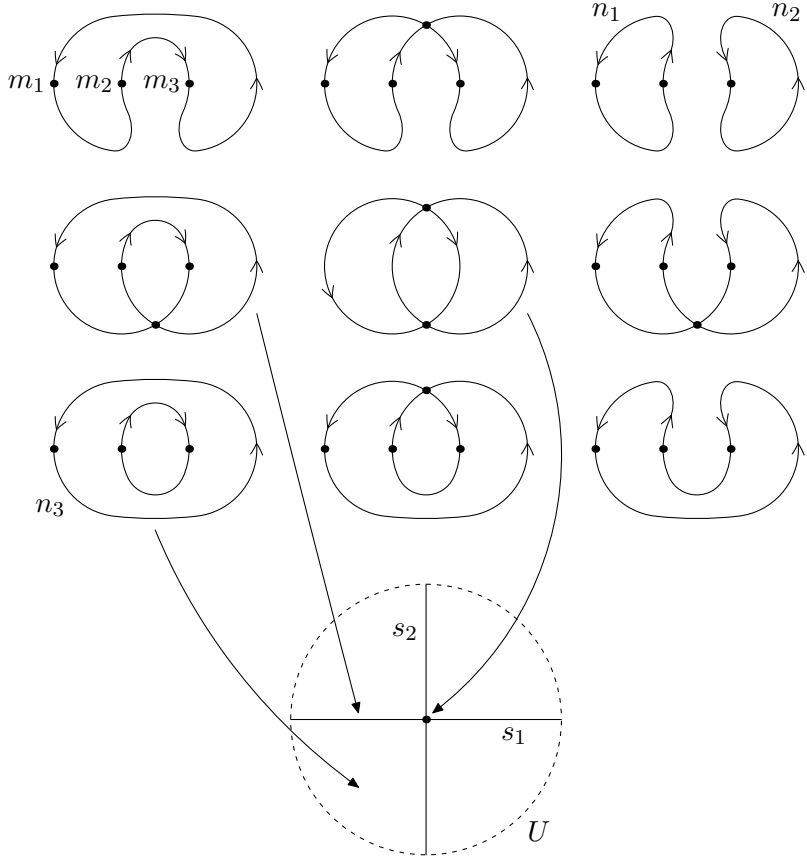}}
  \caption{\label{fig:case2}The same analysis as in
    Figure~\ref{fig:firstsingularity}, in the case when the singular
    fiber is~\usebox{\mygraphics}.  The~$m_i$ are the meridian curves
    of~$H(F)$ and the~$n_i$ are the meridian curves
    of $H(G)$.}
  \label{fig:secondsingularity}
\end{figure}

Let us choose a set of meridian curves $m_1$, $m_2$, $m_3$ for $H(F)$, as
shown in Figure~\ref{fig:secondsingularity}. For $H(G)$ we
pick three curves $n_1$, $n_2$, $n_3$ out of those lying over the two
areas covered by two circles: they form an Heegaard system for
$H(G)$ since they bound disks in it and do not disconnect
$\partial H(G)$. Then the number of intersections between the~$m_i$
and the~$n_i$ is as follows:
\[
\begin{tabular}{Mr|McMcMc}
 m_i \cap n_j  & n_1 & n_2 & n_3 \\ \hlx{vhv}
  m_1 & 1 & 0 & 1 \\
  m_2 & 1 & 0 & 0 \\
  m_3 & 0 & 1 & 0
\end{tabular}
\]
Hence we can reduce the Heegard diagram to the trivial
diagram by eliminating in turn the pairs $(m_3, n_2)$, $(m_2, n_1)$,
and $(m_1,n_3)$.  This shows that in this case as well $H(G)\cup
H(F)=S^3$, with embedded graphs~$G$ (from $\partial_i P''$) and~$F$
(from $\pi^{-1}(x_k)$). In
Figure~\ref{fig:embeddedsecondsingularity} we show how these
graphs and the corresponding meridians are embedded in~$S^3$.

\begin{figure}
  \centerline{\includegraphics{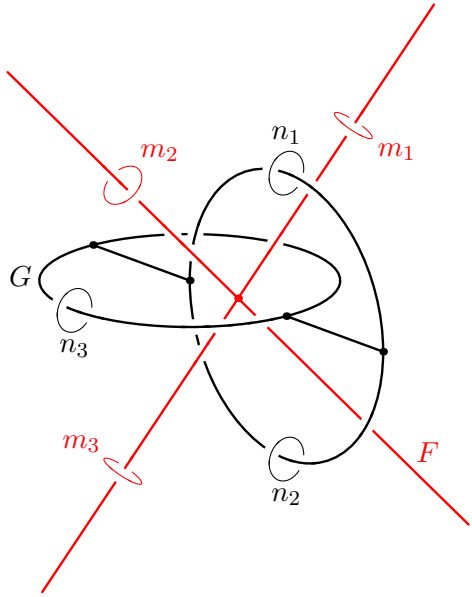}}
  \caption{Two graphs in $S^3$ forming two genus 3 handlebodies
    whose union is the whole space and whose indicated meridians
    intersect as the curves $n_i$ and $m_i$ do. The
    legs of the red graph~$F$ meet in an additional point at $\infty$.}
  \label{fig:embeddedsecondsingularity}
\end{figure}

Now that we know the position of~$G$ in $S^3$, we are left to
construct a simple polyhedron with boundary describing the pair
$(S^3,G)$. We saw in
Example~\ref{examp:secondsingularity} that this graph is represented
by the shadow in Figure~\ref{fig:secondopoliedro}.
Therefore in this case to
extend the construction of the shadow of $M'$ to the singularities of
codimension 2 it is sufficient to use a polyhedron with two
vertices.

\begin{figure}
  \centerline{\includegraphics[width=6.1cm]{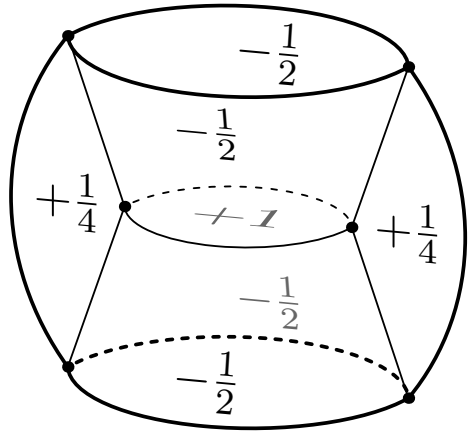}}
  \caption{The polyhedron (containing only two vertices) we use to
    complete the construction of $P'$ near the codimension 2
    singularities of the second type.  The boundary of the
    polyhedron (thicker in the picture) coincides
    with the boundary of the polyhedron $P''$ already constructed away
    from the singularity.}
  \label{fig:secondopoliedro}
\end{figure}

\subsection{Shadow complexity estimate}
\label{sec:estimate-complexity}

In the preceding steps, we constructed a shadow surface~$P'$ for~$M'$,
together with maps $\pi_1:M'\to P'$ and $\pi_2:P'\to \RR^2$ providing
the Stein factorization of $\pi: M' \rightarrow \RR^2$. Let us now
show that the total number of vertices in this shadow is bounded by a
constant times~$t^2$, where $t$~is the number of tetrahedra in the
initial triangulation~$T$.

Given an edge
$e_i$ of~$T$, let $v(e_i)$ be the valence of $e_i$, i.e., the number of
tetrahedra in~$T$ containing $e_i$.  If $e_i$ is contained twice
in a tetrahedron we count it twice. Since each tetrahedron has
$6$~edges, we have
$$\sum_i v(e_i)=6t.$$
Let us now count the total number of segments of singular values
in~$\mr^2$ (i.e., the number of~$f_i$).
In Step 2, while perturbing $\pi$ near an edge $e_i$ in order
to get a stable map, we obtain at most $\frac{v(e_i)}{2}-1$ segments of
critical values.  Thus the total number of~$f_i$ is less than $\sum_i
\frac{v(e_i)}{2}=3t$.  Then, the number of vertices in~$P'$ is
bounded by $2(\#\{f_i\})^2 \le 18t^2$
since each vertex or pair of vertices comes from a crossing.

So the above construction produces a polyhedron~$P'$ with boundary
having a well-controlled number of vertices and admitting a
4-thickening $W_{P'}$ whose boundary contains $M'$.  We
claim that the regions of~$P'$ are disks.  To see this, note that the
regions project by~$\pi_2$ locally homeomorphically onto~$\RR^2$
(since the points in~$P'$ where $\pi_2$ is not locally a homeomorphism
are, by construction, $\Sing(P')$).  Furthermore, by examining
Figures~\ref{fig:firstsingularity} and~\ref{fig:secondsingularity}, we
see that the boundary of the projection of each region turns only to
the left (with the induced orientation from $\RR^2$). Then for each
region of $P'$ apply the following lemma.
\begin{lemma}\label{lem:regionsdiscs}
  Let $R$ be a connected oriented surface with boundary and $\pi:R\to
  \mr^2$ be a local orientation-preserving homeomorphism so that
  $\pi|_{\partial R}$ turns to the left. Then $R$ is a disc and $\pi$
  is an embedding.
\end{lemma}
\begin{prf}
Define the straight arcs in $R$ to be the curves which are
locally projected into straight arcs in $\mr^2$; observe that $\pi$ is
injective on a straight arc.  Let~$p$ be an interior point of $R$.  We
claim that the set of points which can be connected to $p$ by a
straight arc is all of~$R$: this clearly implies the statement.  If
this set were not open then one could find a straight arc connecting
$p$ to another point $q$ whose projection is tangent internally to
$\pi(\partial R)$ but this is ruled out by local convexity of the
image.  On the other hand, because the limit of a sequence of rays is
itself a ray, this set is also closed and so is
all of~$R$.
\end{prf}

Let us now color each edge of $\partial P'$ with one of the colors~$e$
(``external'') and~$f$ (``false''), the ``false'' edges being those
produced during Step~3 of the construction (corresponding to
shrinking circles in~$M'$) and the ``external'' edges being the
remaining ones (around the vertices~$p_i$).  Then~$P'$ is a shadow
for~$M'$, and~$M'$ sits inside the boundary of the thickening~$W_{P'}$
of~$P'$ as the
union of the \emph{horizontal boundary} (i.e., $\pi_1^{-1}(\intr
P')$) and the \emph{vertical boundary} over
$\partial_f P'$ (i.e., $\pi_1^{-1}(\partial_f P')$).  The
components of $\partial M'$ are surfaces~$S_i$ that map
through~$\pi$ to the boundaries of small circles around the~$p_i$.  
By Remark~\ref{rem:genus-boundary}, each spherical component of~$\partial M'$
corresponds to a contractible component of $\partial_e P'$ whose
preimage through~$\pi_1$ in $\partial W_{P'}$ is a 3-ball.

Recall that $M'$ is homeomorphic to $M$ minus a neighborhood of each
non-ideal vertex of~$T$.  Hence, to get a shadow~$P$ for~$M$,
it is sufficient take~$P'$ with the boundary components corresponding
to these non-ideal vertices marked as false.

This concludes the construction of a decorated shadow~$P$ for~$M$ and
an estimate of its shadow complexity.

\subsection{Gleam estimate}
\label{sec:estimate-gleams}

We now provide an upper bound for the absolute value of the gleam of a
region $D$ of $P$.  We remark first of all that if the closure
$\overline{D}$ of $D$ intersects $\partial_f P$, then no gleam has to
be defined on~$D$, as the gleam is defined only on the interior regions
of~$P$.  Note that, by Lemma~\ref{lem:regionsdiscs}, $\overline{D}$
is a closed disc in $P$ with $\partial \overline{D}\subset
\Sing(P)$. Let $v_1,\ldots,v_n$ and $e_1,\ldots,e_n$ be respectively
the vertices and the edges of $\Sing(P)$ contained in $\partial
\overline{D}$ so that $\partial e_i=v_{i+1}\cup v_i$. Let $B$ be a
``shrunk copy of $\overline{D}$'', that is, the complement in
$\overline{D}$ of a small open neighborhood of $\partial \overline{D}$.
Equip $\partial B$ with the decomposition into vertices~$w_i$ and
edges~$f_i$ induced by that of $\partial \overline{D}$.

Now we lift $e_i$ to the strand $e'_i\subset M$ of singular points of
$\pi$ such that $\pi_1(e'_i)=e_i$. This allows us to also choose
lifts~$f'_i$ of $f_i$ by considering the points in $\pi_1^{-1}(f_i)$ which
are nearest to $e'_i$.  Notice that over some~$f_i$, we could have two
choices for~$f'_i$ according to how $D$ is positioned with respect to
the saddle singularity corresponding to~$e_i$.  To get a section of
$\pi_1$ over $\partial B$ we need to choose how to
connect the endpoints of the $f'_i$ with curves in $\pi_1^{-1}
(w_i)$. But $w_i\notin \Sing(P)$ and so
$\pi_1^{-1}(w_i)$ is an oriented curve, so we connect the endpoints of
$f'_i$ and $f'_{i-1}$ through an oriented arc over~$w_i$.  With an
appropriate choice of these arcs over~$w_i$, a section~$s_B$ of
$\pi_1$ over $\partial B$ corresponds to that
used in Proposition~\ref{prop:embpoly} to compute the gleam.

We now construct a section $B'$ of $\pi_1$ over~$B$. To do
this, let $T^{(1)}$ be the $1$-skeleton of~$T$ and for each component
$B_i$ of $B\setminus \pi_1 (T^{(1)})$ choose a face $F_i$ of $T$
so that $B_i\subset \pi_1 (F_i)$ and let $B'_i=\pi_1^{-1}(D_i)\cap
F_i$. One can choose the faces so that they fit coherently on the edges of $T$ projecting inside $B$, because, by construction, these edges correspond to non-critical points. The union of the $B'_i$ is $B'$.

Now compare the section $s_B$ with $\partial B'$ by counting
$\#\{s_B\cap \partial B'\}$. Note that, over an edge $f_i$, $s_B$ and
$\partial B'$ run parallel to each other (because the singular points
of $\pi$ are by construction parallel to the edges of $T$). We claim
that near a vertex $w_i$ they intersect at most once.  Indeed, by
construction $s_B$ is an embedded arc in $\pi_1^{-1}(v_i)$ whose
endpoints are the two endpoints of the lifts $f'_{i}$ and
$f_{i-1}$. On the other side, near $w_i$, $B'$ is formed by a face $F$
of $T$ and at most one of $f'_i$ and $f'_{i-1}$ can be near to an edge
of $F$ so $\partial B'$ can stay near to at most one of $e'_i$ and
$e'_{i-1}$. So, the two sections can intersect over $w_i$ but they can
do it at most once.

Therefore the total gleam of $D$ is bounded above by the number of
vertices in $\partial \overline{D}$ and, since each vertex of $P$
touches exactly $6$ regions, then the gleam weight of $P$ is at most
$6\cdot 18\cdot t^2$.  This concludes the proof of
Theorem~\ref{teo:main}.
\end{prf}

\section{Upper bounds}
\label{sec:upper-bounds}

We can now apply our main tool on constructing shadows,
Theorem~\ref{teo:main}, to give upper bounds on the complexity of
various representations of 3-manifolds in terms of geometric
properties.

\subsection{Triangulating the 4-manifold}
\label{sec:triangulating-4-mfld}
We now turn to the proof of Theorem~\ref{teo:bound}, that 3-manifolds
efficiently bound 4-manifolds.
To actually construct the triangulated 4\hyp manifold, we use the
following:
\begin{lemma}\label{lem:bounding simplices}
Let $(P,g)$ be a boundary\hyp decorated shadow with $n$~vertices whose regions
are disks. There exists a triangulation of
the manifold~$W_P$ obtained by thickening $(P,g)$ (as described by Turaev's
Reconstruction theorem) containing $O(n+|g|)$ simplices. Moreover, the
triangulation can be chosen so that the number of simplices touching
each vertex is bounded above by a constant not depending on $P$ or $g$.
\end{lemma}

\begin{prf}
Choose a triangulation of $P$ such that $\Sing(P)$ is
composed of simplices, the number of simplices depends linearly
(through a constant independent of $P$) on $n$, and the
number of simplices touching each vertex of the triangulation is
bounded above by a constant independent of $P$.  Let $P'$ be the
triangulated polyhedron obtained by deleting a triangle from each
region of $P$ and let $L$ be a 3-dimensional thickening of~$P'$ (which
may be non\hyp orientable). By
gluing prisms on each triangle of $P'$ one can construct a
triangulation of $L$ containing at most $k$ simplices, where
$k=O(n)$. Similarly, the number of simplices touching a vertex can be
linearly bounded from above by a fixed constant not depending on
$(P,g)$.

The 4-dimensional thickening $W_{P'}$ (note that $P'$ has no internal
region and hence no gleam is needed to thicken it) is a fiber bundle
over $L$ with fiber $[-1,1]$. This bundle is the unique bundle
whose total space is orientable, so $W_{P'}$ can be triangulated with
a number of simplices bounded above by $6k$, since $6$ is the minimal
number of simplices needed to triangulate the product of $[-1,1]$ and
a 3-dimensional tetrahedron. Note that the part of the boundary of
$W_{P'}$ which collapses onto the boundary components of
$\partial P'$ that correspond to the triangular punctures is a set of
solid tori $T_i$, all equipped with the same triangulation. In
particular, let $s$ be the number of 3-simplices in the triangulation
of each of these tori.  Applying a Dehn twist to a meridian of a solid
torus~$T_i$,
one obtains a new triangulation of $T_i$ which can be connected by
means of, say, $m$ standard moves of triangulations (called ``Pachner
moves'') to the initial one. Each Pachner move corresponds to gluing a
4-dimensional simplex to $T_i$ along a face and looking at the new
triangulation induced on $T_i$ by the new faces of the simplex. Hence,
in order to perform $g_i$ Dehn twists on $T_i$ it is sufficient to
glue $mg_i$ 4-simplices to $W_{P'}$. Then, gluing a 2-handle on $T_i$
(which can be triangulated with a number of simplices depending only
on $s$) produces a triangulated version of $W_{P}$.

By construction, the number of simplices of dimension~4 in this
triangulation is bounded above by a constant times the number of
simplices in the triangulation of $P'$ plus the sum over all the
regions of $P$ of $s+|g_i|$ where $g_i$ are the gleams and $s$ is the
number of simplices on the solid tori $T_i$.  Moreover, the number of
simplices touching any vertex of the so obtained triangulation can be
controlled from above by a suitable constant which does not depend on
$(P,g)$.
\end{prf}

Theorem \ref{teo:main} and Lemma \ref{lem:bounding simplices} together
prove one main result:
\begin{teo}\label{teo:bound}
  If a 3-manifold $M$ has a triangulation with $t$ tetrahedra, then
  there exists a 4-manifold $W$ such that $\partial W=M$ and $W$~is
  triangulated
  with $O(t^2)$ simplices. Moreover, $W$ has ``bounded geometry'', that
  is, there exists an integer $c$ (not depending on $M$ and $W$) such
  that each vertex of the triangulation of $W$ is contained in less
  than $c$ simplices.
\end{teo} 

An alternate construction gives a simply\hyp connected 4\hyp manifold.

\begin{teo}\label{teo:bound-simply-conn}
A 3-manifold with a triangulation with $k$~tetrahedra is the boundary of 
a simply-connected 4-manifold with $O(k^2)$ 4-simplices.
\end{teo}
\begin{prf}
  Applying the construction underlying the proof of
  Theorem~\ref{teo:main} to the 3\hyp manifold~$M$, we produce a
  decorated shadow~$P$ of~$M$ whose singular set projects
  to~$\mathbb{R}^2$, which contains $O(k^2)$ vertices and whose
  regions are discs.  Note that furthermore that the projection of the
  singular set has $O(k^2)$ self-intersections.  Collapse the false
  boundary of~$P$ and appropriately modify the result to get a genuine
  simple shadow of~$M$ with no more vertices as in
  Remark~\ref{rem:collapse-false-bdry}; the result may not be a
  standard polyhedron, but we can modify it to get a standard
  polyhedron with local moves.  The result is a standard polyhedron,
  which we continue to call~$P$, which still has a map to~$\RR^2$ with
  at most $O(V^2)$ vertices.

  Now apply Proposition~\ref{prop:universallink} to~$P$.  The
  resulting hyperbolic link~$L_P$ is obtained by gluing the pieces as
  in Figure~\ref{fig:blocklink} and performing surgery on some 0\hyp
  framed meridians.  Using the projection to~$\RR^2$, we get a link
  diagram with 1~crossing per vertex of~$P$, at most 3~crossings per
  edge of~$P$, 9~crossings per crossing of the projection of edges
  of~$P$, and 6~crossings per 0\hyp framed meridian: $O(k^2)$
  crossings in all.  $M$~can be obtained by integer surgery on this
  link with coefficients related to the gleams, and so bounded by
  $O(k^2)$.  This surgery diagram therefore gives a simply\hyp
  connected 4\hyp manifold whose boundary is~$M$ with
  complexity~$O(k^2)$.
\end{prf}

Note that the proof of Theorem~\ref{teo:bound-simply-conn} depended on
a reasonably nice projection from the shadow to the plane; in general,
there is no reason to believe that a shadow can be turned into a
simply\hyp connected shadow without an increase in complexity.

\subsection{Upper bounds from geometry}
\label{sec:upper-bounds-geom}
We now turn to hyperbolic geometry, and prove various upper bounds
based on geometry, including an
estimate for the shadow complexity of a hyperbolic 3-manifold by
the square of its Gromov norm.

We first recall an estimate for the number of tetrahedra in a
triangulation from the hyperbolic volume.

\begin{teo}[W.~Thurston]
  \label{teo:spun-triang-vol} There is a positive constant~$C$ so that,
  for every hyperbolic 3\hyp manifold~$M$, there is a
  link~$L$ contained in~$M$ and a partially ideal triangulation of $M
  \setminus L$ with less than $C\cdot \Vol(M)$ tetrahedra.
\end{teo}

This theorem was first used in the proof of the Thurston-J\o rgensen
Theorem~\cite[Theorem 5.11.2]{Thurston82:GeometryTopology}.  We recall the proof here.

\begin{prf}
  Let $V(r)$ be the volume of a ball
  of radius~$r$ in $\HH^3$.

  Let $c$ be the Margulis constant, and consider the thick-thin
  decomposition of~$M$: let~$\eps=c/2$, and let~$M'$ be the
  $\eps$-thick part of~$M$; we will take $L$ to be the core of $M
  \setminus M'$.  Consider a generic maximal $\eps$-net in $M'$, i.e., a
  maximal set of points~$p_i$ such that the distance between any two
  of them is at least $\eps$.  Then the balls of radius~$\eps$ around
  the~$p_i$ cover~$M'$ and the balls of radius~$\eps/2$ around
  the~$p_i$ are embedded in~$M$.  Therefore the number of $p_i$ is
  bounded above by $\Vol(M)/V(\eps/2)$.
  Starting from this
  $\eps$-net we can construct a triangulation~$T$ of a manifold
  homeomorphic to~$M'$ by
  taking the Delaunay triangulation of the~$p_i$; that is, a set
  of~$p_i$ are in a simplex of~$T$ exactly when there is a point~$x$
  in~$M'$ which is equidistant from the chosen~$p_i$ (and no
  farther than $\eps$ from them) and no
  closer to any other~$p_i$.  Since the~$p_i$ were chosen to be
  generic, this is actually a triangulation.

  To estimate the number of tetrahedra in~$T$, let us consider the
  lift~$\tilde T$ of $T$ to $\HH^3$.  Note that two vertices that are
  connected by an edge have distance at most $2\eps$, and so balls of
  radius $\eps/2$ around these vertices fit disjointly in a ball of
  radius $5\eps/2$.  Thus the number
  of vertices connected by an edge to a given vertex in~$\tilde T$ is
  therefore at most $\bigl\lfloor\frac{V(5\eps/2)}{V(\eps/2)}\bigr\rfloor -1$;
  let this number be~$k$.  Then the number of tetrahedra that touch this
  vertex is therefore at most $\binom{k}{3}$ (since in~$\tilde T$ a
  tetrahedron is determined by its set of vertices) and the total
  number of tetrahedra in~$T$ is at most
  \[
  \frac{\Vol(M)}{4V(\eps/2)}\binom{k}{3},
  \]
  which is linear in $\Vol(M)$, as desired.

  To make a partially ideal triangulation of~$M\setminus L$, we need
  to add one more vertex per component of~$L$ and cone from the
  boundary of the triangulation of~$T$ to this new vertex.  This does
  not affect the asymptotic growth of the estimate.
\end{prf}

\begin{teo}\label{teo:gromov-norm-upper}
  There is a universal constant $C>0$ so that a geometric
  3-manifold~$M$, with boundary empty or a union of tori, has shadow
  complexity at most $C\|M\|^2$.
\end{teo}

\begin{prf}
Since Gromov norm is additive under gluing of two
3-manifolds along incompressible boundary tori~\cite{Gromov82:Volume},
by Corollary~\ref{cor:sc-torus-sum} it is sufficient to prove the
theorem for each piece of the JSJ decomposition of~$M$, and by
Proposition~\ref{prop:complexity-zero}
it is sufficient to study the case when
$M$~is hyperbolic.

By applying Theorem~\ref{teo:main} to the triangulation from
Theorem~\ref{teo:spun-triang-vol}, we get a boundary\hyp decorated
shadow~$P_0$ for $M \setminus L$ for some appropriate link~$L$ in~$M$ with
at most $O(\Vol(M)^2)$ vertices.  Collapse the false boundary
of~$P_0$ as in Remark~\ref{rem:collapse-false-bdry} to get a proper
shadow~$P$ for $M\setminus L$; each component
of~$L$ gives a torus boundary component of~$M\setminus L$, which
corresponds to a circle component of~$\partial_e P$.  But now $M$ can
be obtained from $M \setminus L$ by Dehn filling, which by
Proposition~\ref{prop:arbitrarygluing} can be performed at the level
of shadows without adding any vertices.
\end{prf}

\begin{teo}\label{teo:hyperbolic-surgery}
  A finite-volume hyperbolic 3-manifold $M$ with volume~$V$ has a
  rational surgery diagram with $O(V^2)$ crossings.
\end{teo}
\begin{prf}
  The shadow~$P$ constructed in the proof of Theorem~\ref{teo:gromov-norm-upper}
  (before the final Dehn filling) has $O(V^2)$ vertices, and comes
  with a projection to~$\RR^2$ with at most $O(V^2)$ crossings of the
  singular set in the image; as in the proof of
  Theorem~\ref{teo:bound-simply-conn} in the previous subsection, we
  will modify the shadow to give a link diagram.  The only difference
  is that~$P$ has some circular external boundary components
  (corresponding to the link we drilled out) and so cannot be made
  standard; instead, we make all regions either disks or annuli (with
  one external boundary component).  Then we can construct a
  link~$L_P$ as before with at most $O(V^2)$ crossings, and~$M$ is
  surgery on~$L_P$, with rational coefficients on the components
  coming from the external boundary.
\end{prf}

\begin{teo}\label{teo:singularities-upper}
  There exists a constant $C$ such that each hyperbolic 3-manifold $M$
  has a smooth projection in $\mr^2$ with less than $C\|M\|^2$
  crossing singularities.
\end{teo}
\begin{prf}
As in the proof of Theorem~\ref{teo:gromov-norm-upper}, drill out of $M$
some geodesics in order to find a triangulation of a sub-manifold
$M'\subset M$ with a number of tetrahedra bounded above by a constant
times $\Vol(M)$. Applying the construction of the proof of
Theorem~\ref{teo:hyperbolic-surgery} to $M'$, one gets a simple
polyhedron $P$, containing $O(\Vol(M)^2)$ vertices, whose singular set
embeds in $\mr^2$ and such that $M'$ (and hence $M$) is a the Dehn
filling of $S_P$ (recall Proposition~\ref{prop:hyperbolicpieces}). The immersion of 
$\Sing(P)$ in $\mr^2$ produces a smooth projection of
$\pi:S_P\to \mr^2$ with $O(\Vol(M)^2)$ crossing singularities;
moreover, by Lemma \ref{lem:regionsdiscs} the projection through
$\pi$ of each component of $(\partial S_P)$ is a simple curve. 
Hence it is sufficient to show how to
extend $\pi$ to a projection of an arbitrary Dehn filling of $S_P$
without adding any crossing singularities. The idea is to extend $\pi$
through a map whose Stein factorization on the Dehn-filled solid torus
is given by a ``tower'' as that of Figure~\ref{fig:rationalshadow}: in
order not to add any crossing singularities, it is the sufficient to
project the singular set of the added tower to disjoint, nested
circles in $\mr^2$. It is not difficult to check that, even if the
projection of one of these circles intersects the projection of other
components of $\Sing(P)$, no crossing singularity is created since the
fibers in $M$ of the two strands of singular values stay disconnected
around the intersection point.
\end{prf}

\section{Spin boundaries}
\label{sec:spin-boundaries}

In this section we will modify the shadow surfaces constructed in
Section~\ref{sec:shadows-triangulations} to construct a spin-boundary of a
given spin structure on a 3-manifold.  Our goal is the following
theorem.

\begin{teo}\label{teo:bound-spin}
  A 3-manifold with a triangulation with $k$ tetrahedra is the
  boundary of a spin 4-manifold with $O(k^4)$ 4-simplices.
\end{teo}

Before we give the proof, we need some preliminaries about spin
structures on 3- and 4-manifolds given by a shadow.

\subsection{Spin structures from shadows}
\label{sec:spin-structures}

We first identify the Stiefel-Whitney class $w_2(W) \in H^2(W; \ZZ/2)$
for a 4-manifold $W$ with shadow $P$.  In a 4\hyp manifold, the
evaluation of $w_2$ on a closed surface $S$ is the reduction modulo 2
of the self-intersection number of $S$.  Self-intersection numbers for
a general integral homology class may be computed with the following
proposition.

\begin{prop}[Turaev, \cite{Turaev91:TopologyShadows}]
  Given a shadow $P$ for a 4-manifold $W$ with oriented regions $f_i$
  and gleams $g(f_i)$, and $S \in H_2(W; \ZZ)$ given by the chain $S =
  \sum_i a_i f_i$.  The self\hyp intersection number of $S$ is
  \[
  S\cdot S = \sum_i a_i^2 g(f_i).
  \]
\end{prop}

For an element of $H_2(W; \ZZ/2)$, this simplifies: a $\ZZ/2$ homology
class is a formal sum of regions of the shadow surface, with an even
number of regions (0 or 2) around each edge; that is, it is a union of
regions of the shadow surface which form a closed surface without
singularities.  Since $a_i^2 = a_i$ in $\ZZ/2$, we have:

\begin{prop}\label{prop:self-intersection}
  The self-intersection number of $S \in H_2(W; \ZZ/2)$ is the sum of
  gleams on the regions that appear in $S$.
\end{prop}

Note that in both cases it is not immediately obvious that the sum is
an integer.  If the gleams were all integers, this proposition would
say that the Stiefel-Whitney class is the reduction modulo 2 of the
cocycle given by the gleams.  In general, the following holds:
\begin{lemma}
  \label{lem:stiefel-whitney-integer}
  The $\ZZ/2$-cochain represented by the $\ZZ/2$-gleam is a
  coboundary.  Therefore the gleam cocycle $g$ is cobordant (over
  $\QQ$) to an integer-valued cocycle $g'$, for each cycle $z=\sum_i a_i
  f_i\in H_2(P;\ZZ)$ the sum $\sum_i g(f_i) a_i$ is an integer, and the
  Stiefel-Whitney class is the reduction modulo 2 of $g'$.
\end{lemma}
\begin{prf}
We will explicitly construct a 0-cochain with coefficients in $\ZZ/2$
whose coboundary is the $\ZZ/2$-gleam cochain. First define an
orientation on a vertex $v$ of a simple polyhedron to be a numbering
(from 0 to 3) of the 4 edges of the singular set touching the
vertex. An orientation on a vertex allows one to name the six regions
touching the vertex as $R_{\{i,j\}}$ where $i\neq j$ are in
$\{0,1,2,3\}$ so that the edge $k$ is touched by the regions
$R_{\{k,j\}}$ with $k\neq j$.  Construct a cyclic
ordering of the three regions touching each edge by considering the
order induced by the $j$'s if~$k$ is odd and the reverse if~$k$ is even.

Now fix arbitrarily an orientation around each vertex of $P$ and
consider the 1-cochain $c$ whose value on an edge $e\in \Sing(P)$ is 1
if the two cyclic orderings induced on the three regions touching $e$
by the two vertices touched by $e$ are the same and 0 otherwise.  For
an edge $e$ which is a circle, $c(e)$ is 0 if there is a consistent
cyclic orientation on the regions touching $e$ and 1 otherwise.  We
claim that $\delta c$ is the $\ZZ/2$-gleam cochain.  For simplicity,
we prove it on a region $R$ whose boundary passes at most once on each
edge of $\Sing(P)$. In that case, the $\ZZ/2$-gleam is 1 iff the
regular neighborhood of $\partial R$ in $P\setminus R$ collapses on an
odd number of M\"obius strips.  For each component $\alpha$ of
$\partial R$, let us fix an orientation of $\alpha$ and a base point
$p_{\alpha}$ contained in an edge $e_{\alpha}$ of $\Sing(P)$ and lying
in a neighborhood of a vertex $v_{\alpha}$. To control the topology of
the neighborhood of~$\alpha$ it is sufficient to count (mod~$2$) 
how many times, while running along~$\alpha$, 
two consecutive vertices are connected 
by an orientation preserving gluing. This number is $\langle\delta
c,\alpha\rangle$.  This proves the first statement.

For the second part, define $g'$ by
\[
  g'(f_i) = g(f_i) + \frac{1}{2}\delta c(f_i),
\]
where $c$ is the $\ZZ$-valued 1-cochain defined as above, considering
0 and 1 as elements of $\ZZ$ rather than $\ZZ/2$. By the above
results, $g'(f_i)$ are integers.  Since $g$ and $g'$ are cobordant,
for a closed surface $S = \sum a_if_i$, we have $\sum_ig(f_i)a_i =
\sum_ig'(f_i)a_i$.  This second sum is obviously an integer and so, by
Proposition~\ref{prop:self-intersection}, $g'$ represents the
Stiefel-Whitney class.
\end{prf}

For a given 3-manifold~$M$, we are interested not just in finding a
spin 4-manifold~$W$ with $\partial W = M$, but in finding one where a given
spin structure~$s$ on~$M$ extends to a spin structure on~$W$.  The
obstruction to extending a given spin structure is a relative
Stiefel-Whitney class $w_2(W, s) \in H^2(W, \partial W; \ZZ/2)$.  Since $H^2(W,
\partial W; \ZZ/2) \cong H_2(W; \ZZ/2)$ by Poincar\'e duality, $w_2(W,s)$ is
given by a subsurface~$F$ of $P$, possibly disconnected or not
orientable.  The image of $w_2(W,s)$ under the
natural map from $H^2(W, \partial W)$ to $H^2(W)$ must be the original
Stiefel-Whitney class of $W$.  The corresponding map from $H_2(W)$ to
$H^2(W)$ is given on a homology cycle by the cup product, so for any
other subsurface~$S$ of~$P$, we must have $S\cdot F = S\cdot S$.  A
surface~$F$ satisfying this property is called a \emph{characteristic
  surface}.  Every characteristic surface appears as the obstruction
to extending a spin structure on the boundary.

\subsection{Constructing efficient spin fillings}
\label{sec:spin-fillings}
With these preliminaries in hand, we prove the following result which
allows us to convert an arbitrary shadow to one that spin-bounds a
given spin structure.  Together with
Theorem~\ref{teo:main} and Lemma~\ref{lem:bounding simplices} this implies
Theorem~\ref{teo:bound-spin}.
\begin{teo}\label{teo:spin-estimate}
  Let $M$ be a 3-manifold equipped with a spin-structure~$s$, let $W$
  be a 4-manifold with $\partial W=M$ and let $(P,g)$ be a shadow
  of~$W$ with $k$~vertices and gleam weight~$|g|$. Then there exists a
  shadow~$P'$ of~$M$ containing $O(k^2)$ vertices, with gleam weight
  $|g| + O(k)$, and whose thickening is a 4-manifold~$W'$ admitting a
  spin-structure whose restriction to $\partial W=M$ is~$s$.
\end{teo}

\begin{prf}
Let us first sketch the strategy of the proof.  As outlined above, the
Poincar\'e
dual of $w_2(W,s)$ is a 2-cycle represented by a surface $F$ in $P$.
Let $N(F)$ be a regular neighborhood of $F$ in $W$ (diffeomorphic to a
disk bundle over~$F$); by construction, there exists a spin-structure
on $W\setminus N(F)$ extending $s$. Note that $P\cap (W\setminus
N(F))$ is a polyhedron with boundary, embedded in $W\setminus N(F)$,
with boundary in $\partial N(F)$.  Our goal is to replace $N(F)$ with
a suitable 4-manifold $W_F$ equipped with a shadow with boundary, so
that the spin-structure of $W\setminus N(F)$ extends on $W_F$ and that
the gluing of the shadows matches up the boundaries.  We will do this
by dividing~$N(F)$ into a number of pieces, of a finite number of
topological types, each of which can be replaced by a standard model.
There are some special regions to incorporate the gleams.  In all this
process we will have to control the complexity of the resulting final
shadow.

The proof is divided into four steps:
\begin{enumerate}
\item Estimate the complexity of the shadow $P\cap N(F)$;
\item Find an efficient set of separating curves in $F$ decomposing it
  into M\"obius strips,  punctured spheres and tori;
\item Perform surgery along these curves producing $W'$ and its shadow
  $P'$, reducing to the case when
  $F$ is a union of $\mathbb{RP}^2$, $S^2$ and $T^2$;
\item Solve the cases for $F=S^2$, $\mathbb{RP}^2$, and $T^2$.
\end{enumerate}
\step 1
The polyhedron with boundary $P\cap N(F)$ is properly embedded in $N(F)$
and is the mapping cylinder of the projection of the trivalent graph
$G=P\cap \partial N(F)$ into $F$. The set $\Sing(P)\cap F$ coincides with
the graph formed by the projection of $G$ in $F$ and its vertices
correspond to distinct vertices of $P$ and have valency either 3
(corresponding to vertices of $G$) or 4 (corresponding to transverse
self intersections of the projection of $G$ in $F$). In
particular, it contains $O(k)$ vertices and edges. $F$ is a union of
regions of $P$ and so comes equipped with the gleams induced by the
inclusion; in particular, the gleam weight of $P\cap N(F)$ is bounded by
$|g|$.

\step 2
By the preceding step, the graph $E=\Sing(P)\cap F$ has $O(k)$
vertices with valency 3 or 4 and splits $F$ into regions with gleam
whose total weight is no more than~$|g|$.  We will now define a
$0\to2$-move on a pair of edges $e_0$ and $e_1$ of $E$. Let $R_0$ and
$R_1$ be the regions of $P\setminus F$ such that $\overline{R_i}\cap
F=e_i,\ i=0,1$. A $0\to 2$-move on the pair $\{e_0,e_1\}$ is the
sliding $R_0$ over a small disk contained in $R_1$ (the construction
is symmetrical) as shown in Figure \ref{fig:0-2forspin}. In the lower
part of the figure we exhibit the modification induced on $E$ by the
move.
\begin{figure}[htbp]
 \centerline{\includegraphics[width=8.4cm]{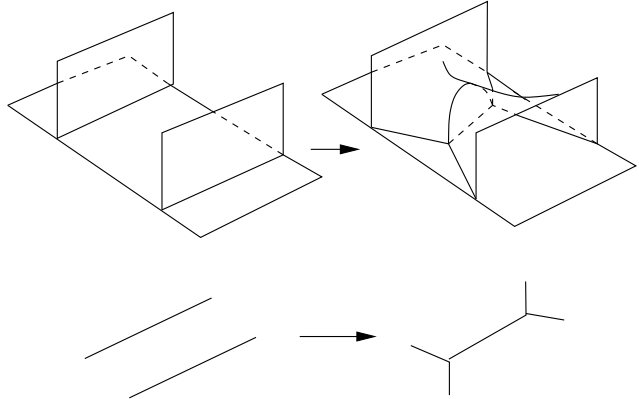}}
  \caption{The $0\to 2$-move and its effect on the graph $\Sing(P)\cap F$.}
  \label{fig:0-2forspin}
\end{figure}
A $0\to 2$ move does not modify $W$ nor $F$ but changes $P$ by
adding 2 vertices.

\begin{lemma}\label{lem:induction}
Let $S$ be a connected surface and $H\subset S$ be a graph containing
$n$~vertices each having valency~3 or~4. After applying
$O((|\chi(S)|+n)^2)$ $0\to 2$-moves on~$H$, it is possible to find a
set~$C$ of closed curves, with each component of~$C$ a separating
curve, cutting the pair $(F,H)$ into the following blocks:
\begin{enumerate}
\item Pairs $(D^2,H)$, where $H$ is a properly embedded graph in $D^2$
  with one vertex of valency 3 or 4.
\item Pairs $(A,H_A)$ where $A$ is an annulus and $H_A$ is a trivalent
  graph in $A$;
\item Pairs $(P,H_P)$ where $P$ is a thrice punctured sphere and $H_P$
  is a (possibly empty) set of disjoint, non-parallel, essential arcs
  in $P$;
\item Pairs $(M,H_M)$ where $M$ is a M\"obius band and $H_M$ is either
  empty or an essential arc in $M$;
\item Pairs $(T,H_T)$ where $T$ is a punctured torus and $H_T$ is a
  (possibly empty) set of disjoint, non-parallel essential arcs in~$T$.
\end{enumerate} 

Moreover the set $C$ can be chosen to intersect $H$ in
$O((|\chi(S)|+n)^2)$ points.

\end{lemma}
\begin{prf}
If $H$ is empty, the lemma is trivial.  If~$H$ has no vertices, we can
create 2~vertices by a $0 \to 2$ move.

Encircle each vertex of $H$ by a closed curve: this set
of $n$ curves intersects~$H$ at most $4n$ times and decomposes~$S$
into $n$ blocks of the first type and a surface~$S'$ whose Euler
characteristic is $\chi(S)-n$. If $H'=H\cap S'$ contains parallel
edges we apply $O(n)$ $0\to 2$-moves in order to replace each set of
parallel edges by a single edge branching only in the neighborhood of
$\partial S'$. Then we add one curve per component of $\partial S'$
in order to enclose all these trivalent vertices in annular
regions. This creates at most $n$~annular regions of type~2 in the
above list; the set of curves considered until now intersects~$H$ at
most $8n$ times.
Let $S''\subset S'$ be the remaining surface equipped with the graph
$H''=H'\cap S''$.
We will prove the statement of the lemma for the pair $(S'',H'')$ by
arguing by induction on $|H_1(S'';\mz_2)|$.
If $\chi(S'')=0$ or if $S''$ is orientable and $\chi(S'')=-1$ we are
done.  Otherwise the result follows by induction by applying
Lemma~\ref{lem:spin-inductive-step} below to~$S''$.  Note that
$\chi(S'') = \chi(S) - n$.
\end{prf}

\begin{lemma}\label{lem:spin-inductive-step}
  Let $S$ be a connected surface with boundary, with
  $\chi(S) < 0$ and not a thrice-punctured sphere, equipped with a
  set~$H$ of essential, pairwise non-parallel arcs.
  Then it is possible to apply $O(-\chi(S))$ $0\to 2$-moves to~$H$ and
  find a set of curves in~$S$ intersecting~$H$ at most $-18\chi(S)$
  times and cutting~$S$ into annuli of type~2 in the above list and
  two surfaces~$S_1$ and~$S_2$ such that
  $|H_1(S_i;\mz_2)|<|H_1(S;\mz_2)|$ for $i=1,2$.
\end{lemma}

\begin{prf}
  We first find an essential separating curve~$c_0$ cutting~$S$ into~$S_1$
  and~$S_2$ and intersecting~$H$ at most $-4\chi(S)$ times. Then
  we show how to apply $0\to 2$-moves and add boundary parallel
  curves in order to prove the rest of the claim.

  Complete~$H$ to a maximal set of essential, pairwise non-parallel
  arcs, that is, an ideal triangulation of $S$. This can be achieved
  by adding $O(\chi(S))$ arcs: indeed, the cardinality of $H$ after
  this becomes exactly $-3\chi(S)$.  We now prove that~$c_0$ can be
  chosen to intersect $H$ at most $-12\chi(S)$
  times. By cutting $S$ along $1-\chi(S)$ arcs of $H$ and shrinking
  the original boundary components to points, we can reduce to
  considering a polygon with $2-2\chi(S)$ sides whose sides are
  identified in pairs.  The other arcs of~$H$ are contained among the
  diagonals of the polygon.

  If two edges of the polygon are identified in~$S$ through an
  orientation-preserving homeomorphism, then the straight curve
  connecting their midpoints is a simple closed curve in~$S$ whose
  regular neighborhood is a M\"obius strip whose boundary is~$c_0$: it
  intersects each edge of~$H$ at most two times and it cuts a M\"obius
  strip out of~$S$.  From now on we will assume all edges are
  identified by an orientation-reversing homeomorphism, i.e., $S$~is
  orientable.

  If there is a pair $a, a'$ identified in~$S$ such that the straight
  arc~$\alpha$ connecting their midpoints is not separating in~$S$,
  then we can find another pair $b, b'$ identified in~$S$ such that
  the straight arc~$\beta$ connecting their midpoints
  intersects~$\alpha$ once. Then the curve~$c_0$ formed by the
  boundary of the regular neighborhood in~$S$ of $\alpha\cup \beta$
  intersects each edge of $H$ at most four times and cuts a punctured
  torus out of~$S$.

  If there is a pair of edges $a, a'$ so that the straight
  arc~$\alpha$ connecting their midpoints is an essential,
  disconnecting curve in~$S$, then we can take~$c_0$ to be~$\alpha$.
  In this case it intersects each edge of~$H$ at most once.

  If none of the above cases hold, each edge is paired with a
  neighboring edge.  In this case, pick four neighboring edges $a, a',
  b, b'$, with $a, a'$ and $b, b'$ paired, and take $c_0$ to be the
  union of the straight arcs connecting the midpoints of $a'$ and $b$
  and $a$ and $b'$.  In this case it intersects each edge of~$H$ at
  most twice.

  In all cases the curve~$c_0$ intersects $H$ at most $4\# H\le
  -12\chi(S)$ times and cuts $S$ into the union of two surfaces
  $S_1\cup S_2$ such that $|H_1(S_i;\mz_2)|< |H_1(S'';\mz_2)|,\
  i=1,2$. We still need~$H_i$ to be composed only of non-parallel,
  essential edges in~$S_i$. To fulfill this condition, it is
  sufficient to apply to each pair of parallel edges of~$H_i$ a
  $0\to 2$-move in order to produce trivalent graphs in $S_i$ each
  having its vertices near~$\partial S_i$. Then we add to~$c_0$ a copy
  of each of the components of~$\partial S_i$, creating at most $-6\chi(S)$
  intersections, in order to enclose these trivalent
  vertices in annuli.  (In order to count the number of new intersections,
  note that after joining parallel arcs, in~$S_i$ there are at most
  $-3\chi(S_i)$ essential arcs, each of which touches $\partial S_i$
  twice, and $\chi(S_1) + \chi(S_2) = \chi(S)$.)
  The remaining blocks are surfaces $S'_i,\ i=1,2$, homeomorphic to
  $S_i$ and equipped with sets of disjoint essential arcs $H'_i,\
  i=1,2$.
\end{prf}

Applying Lemma~\ref{lem:induction} to each connected component of
$(F,E)$ and noting that the sum of $|\chi(F_i)|$ over the
components~$F_i$ of~$F$ is $O(k)$, we conclude that, after applying
$O(k^2)$ $0\to 2$-moves to the edges of~$E$ (which adds $O(k^2)$
vertices to~$P$), it is possible to find a set of disjoint, separating
curves cutting~$F$ into disks, annuli, once punctured $T^2$ and
$\mathbb{RP}^2$ and thrice punctured spheres, whose total number of
intersections with~$E$ is~$O(k^2)$.

Moreover, by adding one simple curve bounding a disk enclosing as much
of the gleam as possible per region, we can suppose that the gleam of
other regions is~$0$ or~$\pm \frac12$, depending on the $\ZZ/2$
gleam.  We do not add these curves to the regions with zero gleam.
These circles will be called \emph{gleam circles}.

\step 3
We will now perform 4-dimensional surgeries inside
$N(F)$ in order to replace a regular neighborhood of each curve of the
family~$C$ found in the preceding step with the regular neighborhood of a
self-intersection~0 sphere. Topologically, this corresponds to
replacing $S^1\times D^3$ with $S^2\times D^2$ (whose boundary is
$S^2\times S^1$ in
both cases). On the level of polyhedra, each move replaces the regular
neighborhood of a separating curve in~$F$ with a polyhedron collapsing
on a sphere whose total gleam is 0 as shown in
Figure~\ref{fig:surgeryonpoly}. 
\begin{figure}
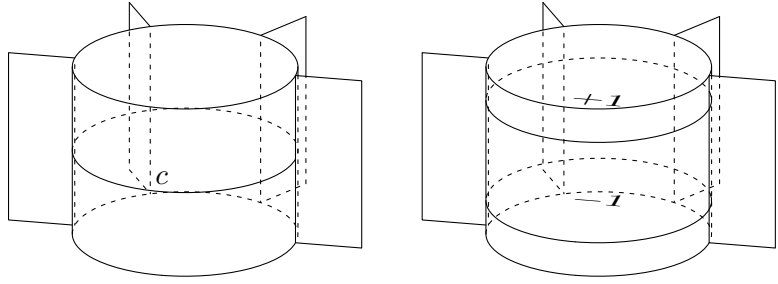

  \[
  \mfig{draws/shadow.40}\qquad\mfig{draws/shadow.41}
  \]
  \caption{On the left we draw a typical regular neighborhood of a
    separating curve $c\subset F$: in the drawing $F$ is the vertical
    cylinder and the lateral regions correspond to the regions of $P\setminus
    F$ whose intersection with $F$ form $E$. 
  On the right we show the result of a surgery on $c$: two new disks
  are glued to $F$ along curves parallel to $c$; their gleams are
  respectively 1 and $-1$. }
  \label{fig:surgeryonpoly}
\end{figure} 
After modifying $N(F)$ and $P\cap N(F)$ in this way, we get a new
4-manifold $N'(F)$ equipped with a shadow $P'$ with the same boundary
as $N(F)$ and $P\cap N(F)$.  Since the set of curves
over which we performed the surgeries intersects $\Sing(P)\cap F$ only
$O(k^2)$ times, the number of vertices of~$P'$ is $O(k^2)$.  Moreover
the gleam weight of $P'$ is increased by
$O(k)$ because of the addition of the two disks of gleam $\pm 1$ for
each curve over which we performed a surgery. 
\begin{lemma}\label{lem:surgery-poly}
  After surgery on a separating curve~$c$ of the characteristic surface as
  above, $w_2(N(P),s)$ is represented by the
  sub-polyhedron~$F'$ of $P'$ obtained by cutting $F$ along all the curves
  over which we performed the surgeries and capping the resulting
  boundaries with the $\pm 1$-gleam disks added during the surgeries.
\end{lemma}
\begin{prf}
  The characteristic surface agrees with $F$ away from the surgery,
  and so can either be the original surface $F$, including the annulus
  between the two disks, or the surface $F'$, not including the
  annulus but including the two disks.  Let $F_1$ be one of the two
  components of $F\setminus c$, and let $g_1$ be its total gleam (which is
  necessarily an integer).  Consider the surface $F_1'$ obtained by
  attaching to $F_1$ the disk with gleam $+1$ after surgery.  The
  self-intersection number of $F_1'$ is $F_1'\cdot F_1' = g_1+1$, while
  $F\cdot F_1' = g_1$, so $F$ itself is not characteristic and the new
  characteristic surface must be $F'$.
\end{prf}

By Lemma~\ref{lem:surgery-poly}, after surgery on a complete set of
separating curves the characteristic surface is a surface $F'$
formed by a disjoint union of surfaces~$S_i$, $i=1,\ldots m$ each of
which is a $S^2$, $\mathbb{RP}^2$ or $T^2$ embedded in $P'$, which we
will call~$P$.

Note that, by the construction in Step~2, almost all of the regions
of~$F'$ are equipped with gleam at most $\pm \frac12$.  The regions
possibly having non-zero gleam are those coming from the disks added while surgering
or from the disks bounded by the gleam circles.
Let us call the spheres that result from these
last disks \emph{gleam spheres}.  Therefore each~$S_i$ which is not a
gleam sphere is either a $S^2$ with total gleam in $[-3,3]$ or a $T^2$
or $\mrp^2$ with total gleam in $[-2,2]$.

\step{4} We have now reduced the problem to the case when $F$ is a
disjoint union of gleam spheres and surfaces of a finite number of
types with bounded gleams.  Note that, by the construction of Step~2,
on each~$S_i$ which is not a gleam sphere the graph $E_i=\Sing(P)\cap
S_i$, in addition to the gleam circles (each of which now bounds a
disk of gleam $\pm 1$), is one of the following graphs:
\begin{enumerate}
\item $S_i$ a sphere coming from a vertex of the original
  singular graph: $E_i$~is a circle and the cone from its
  center to 3 or 4 points. $S_i$ is split into 4 or 5 disks,
  one of which has gleam $\pm 1$ and the other ones having gleam~0
  or~$\pm\frac12$.
\item \label{item:annuli-sphere}
  $S_i$ a sphere coming from an annulus:
  $E_i$ is two circles bounding $\pm 1$ gleam disks and a
  trivalent graph connecting them whose complement is made of gleam~0
  regions without gleam circles.
\item $S_i$ a sphere coming from a thrice-punctured sphere: $E_i$ is
  three circles bounding three $\pm 1$-gleam disks in $S_i$, connected
  by a subset of the three essential arcs
  in the complement of the disks, cutting it into gleam~0 regions.
\item $S_i$ is a torus or a $\mrp^2$: $E_i$ is formed
  by a circle bounding a $\pm 1$-gleam disk and a set of disjoint, essential
  arcs in the complement of this disk, cutting it into gleam~0
  regions.
\end{enumerate}
\begin{figure}
 \centerline{\includegraphics[width=8.4cm]{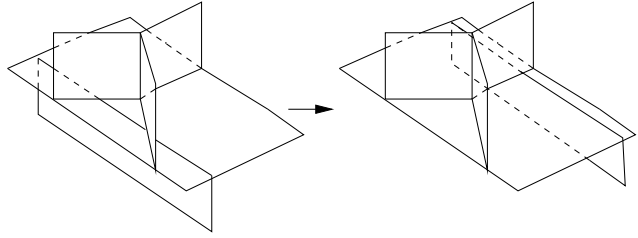}}
 \caption{The $3\to 2$-move decreases the number of vertices of $P$
   without changing its thickening. }
  \label{fig:2to3forspin}
\end{figure} 

Now applying some $3\to 2$-moves to~$P$ as shown in
Figure~\ref{fig:2to3forspin} (which decreases the
number of vertices), we can push all the trivalent vertices of~$E_i$
out of the spheres of type~\ref{item:annuli-sphere} above, so to
reduce them to spheres in which $E_i$ is made of two circles
bounding $\pm 1$-gleam disks, connected by some number of parallel
arcs.  If there is more than one parallel arc, we can reduce the
number of arcs by performing a $0 \to 2$ move followed by two $3 \to
2$ moves, not changing the total number of vertices.  So we may assume
in this case that $E_i$ is two circles, possibly connected by one arc.

Thus, each pair $(S_i, E_i)$ which is not a gleam sphere is one of a
finite number of cases.  Since every spin 3-manifold spin-bounds a
4-handlebody and every 4-handlebody admits a shadow, for each of these
cases we can choose a spin-filling equipped with a shadow.  Therefore
we replace the regular neighborhood of each such~$S_i$ in $W$ with the
suitable model. It is clear that this changes the complexity of $P$ by
a fixed finite amount, which can be estimated by explicitly choosing
the above models for the spin fillings.

We are left with the gleam spheres, each of which is a sphere split by
a simple closed curve into a $\pm 1$-gleam disk and a $g$-gleam disk
(for suitable $g$).  We can split these last spheres by a surgery
along a curve inside the $g$-gleam disk bounding a disk of gleam
$g-\sign(g)$. By surgering over this curve and iterating we reduce to
a union of $O(|g|)$ spheres each being composed of two $\pm 1$-gleam
disks and one $\pm 1$-gleam annulus, which, again, can be replaced by
a fixed spin-filling chosen once and for all.

This concludes Step~4.
\end{prf}

\bibliographystyle{hamsplain}
\bibliography{topo,quantum}

\end{document}

%% file: macros.tex
\def\mathcenter#1{%
  \vcenter{\hbox{#1}}%
}

\def\mfig#1{
        \mathcenter{\includegraphics{#1}}
}

\def\mfigb#1{
        \mathcenter{\includegraphics[trim=-1 -1 -1 -1]{#1}}
}

\newcommand{\case}[2]{%
  \smallskip%
  \noindent\textbf{Case #1.} #2\par\nobreak%
}

\newtheoremstyle{citing}% name
  {}%      Space above, empty = `usual value'
  {}%      Space below
  {\itshape}% Body font
  {}%         Indent amount (empty = no indent, \parindent = para indent)
  {\bfseries}% Thm head font
  {.}%        Punctuation after thm head
  {.5em}%     Space after thm head: " " = normal interword space;
  {\thmnote{#3}}% Thm head spec

\theoremstyle{citing}
\newtheorem*{citingthm}{} % All text appears in the note
\theoremstyle{plain}

%% file: defs.tex
\newtheorem{teo}{Theorem}[section]
\newtheorem{theorem}[teo]{Theorem}
\newtheorem{lemma}[teo]{Lemma}
\newtheorem{prop}[teo]{Proposition}
\newtheorem{cor}[teo]{Corollary}

\theoremstyle{definition}
\newtheorem{defi}[teo]{Definition}
\newtheorem{definition}[teo]{Definition}
\newtheorem{rem}[teo]{Remark}
\newtheorem{example}[teo]{Example}

\newtheorem{question}[teo]{Question}

\theoremstyle{remark}

\newenvironment{prf}[1][Proof.]{%
  \begin{proof}[{#1}]%
  \edef\qedsymbol{\noexpand\fbox{{\noexpand\Small\theteo}}}%
}{\end{proof}}


\newcommand{\RR}{\mathbb R}

\newcommand{\ZZ}{\mathbb Z}
\newcommand{\QQ}{\mathbb Q}
\newcommand{\HH}{\mathbb H}
\newcommand{\mr}{\mathbb{R}}
\newcommand{\mz}{\mathbb{Z}}
\newcommand{\mrp}{\mathbb{RP}}
\newcommand{\tet}{\mathrm{tet}}
\newcommand{\oct}{\mathrm{oct}}

\newcommand{\eps}{\varepsilon}

\DeclareMathOperator{\sign}{sign}

\DeclareMathOperator{\vol}{vol}
\DeclareMathOperator{\Sing}{Sing}
\DeclareMathOperator{\intr}{int}

\DeclareMathOperator{\scomp}{sc}
\DeclareMathOperator{\Vol}{Vol}

\newcommand{\connect}{\mathbin \#}

\newcommand{\step}[1]{\medskip\par\noindent\textbf{Step #1.}}